\documentclass[12pt,sort,compress]{elsarticle}

\usepackage{amsmath}
\usepackage{amssymb}
\usepackage{amsfonts}
\usepackage{amsthm}
\usepackage{mathrsfs}
\usepackage{bm}
\usepackage{nicefrac}
\usepackage{siunitx}
\usepackage{xfrac}

\usepackage{graphicx}
\usepackage{epsfig}
\usepackage{float}
\usepackage{textcomp}
\usepackage{graphbox}
\usepackage{wrapfig}
\usepackage{subcaption}
\usepackage{caption}

\usepackage[top=1in,bottom=1in,left=1in,right=1in]{geometry}
\usepackage{color}
\usepackage{xcolor}
\usepackage{url}
\usepackage{multicol}
\usepackage{parskip}
\usepackage{titlesec}
\usepackage{textcomp}
\usepackage{textgreek}
\usepackage{times}
\usepackage{enumitem}
\usepackage{array}
\usepackage{arydshln}
\usepackage{booktabs}
\usepackage{lipsum}
\usepackage{gensymb}
\usepackage{nth}
\usepackage{setspace}

\usepackage{multirow}
\usepackage{tabulary}

\usepackage{algorithm}
\usepackage{algorithmic}
\usepackage{listings}
\usepackage{manyfoot}

\usepackage{tikz}
\usepackage{pgfplots}
\usepackage{pgfplotstable}
\pgfplotsset{compat=newest,scaled y ticks=true}
\usepgfplotslibrary{units}
\usetikzlibrary{spy}
\usetikzlibrary{backgrounds}
\usetikzlibrary{shapes,arrows,positioning,intersections,quotes}
\usetikzlibrary{matrix}

\usepackage[pdfborder={0 0 0},colorlinks,allcolors=blue]{hyperref}

\graphicspath{{Figures/},{New_Figures/},{figs/}}
\DeclareGraphicsExtensions{.pdf,.png,.jpg}

\theoremstyle{definition}

\allowdisplaybreaks

\begin{document}
\begin{frontmatter}


\title{\textsc{ValveFit}: An analysis-suitable B-spline-based surface fitting framework for patient-specific modeling of tricuspid valves}

\author[1]{Ajith~Moola}\ead{arma@iastate.edu}

\author[1]{Ashton~M.~Corpuz}\ead{amcorpuz@iastate.edu}

\author[2]{Michael~J.~Burkhart}\ead{michael.j.burkhart-1@ou.edu}

\author[2]{Colton~J.~Ross}\ead{colejr98@gmail.com}

\author[3]{Arshid~Mir}\ead{arshid-mir@ouhsc.edu}

\author[4]{Harold~M.~Burkhart}\ead{harold-burkhart@ouhsc.edu}

\author[5]{Chung-Hao~Lee}\ead{chunghao.lee@ucr.edu}

\author[1]{Ming-Chen~Hsu}\ead{jmchsu@iastate.edu}

\author[1]{Aishwarya~Pawar\corref{cor1}}\ead{arpawar@iastate.edu}

\address[1]{Department of Mechanical Engineering, Iowa State University, Black Engineering Building, 2529 Union Drive, \\Ames, Iowa 50011, USA}

\address[2]{School of Aerospace and Mechanical Engineering, The University of Oklahoma, Felgar Hall, 865 Asp Avenue, \\Norman, Oklahoma 73019, USA}

\address[3]{Department of Pediatrics, The University of Oklahoma Health Sciences Center, 1200 N. Children's Ave., Suite 2F, \\Oklahoma City, Oklahoma 73104, USA}

\address[4]{Department of Surgery, The University of Oklahoma Health Sciences Center, 800 Stanton L. Young Blvd., \\Oklahoma City, Oklahoma 73104, USA}

\address[5]{Department of Bioengineering, University of California, Riverside, 900 University Ave., \\Riverside, California 92507, USA}

\cortext[cor1]{Corresponding author}

\journal{Computer Methods in Applied Mechanics and Engineering}

\begin{abstract}
Patient-specific computational modeling of the tricuspid valve (TV) is vital for the clinical assessment of heart valve diseases. However, this process is hindered by limitations inherent in the medical image data, such as noise and sparsity, as well as by complex valve dynamics. We present \textsc{ValveFit}, a novel differentiable B-spline surface fitting framework that enables rapid reconstruction of smooth, analysis-suitable geometry from point clouds obtained via medical image segmentation. We start with an idealized TV B-spline template surface and optimize its control point positions to fit segmented point clouds via an innovative loss function, balancing shape fidelity and mesh regularization. Novel regularization terms are introduced to ensure that the surface remains smooth, regular, and intersection-free during large deformations. We demonstrate the robustness and validate the accuracy of the framework by first applying it to simulation-derived point clouds that serve as the ground truth. We further show its robustness across different point cloud densities and noise levels. Finally, we demonstrate the performance of the framework toward fitting point clouds obtained from real patients at different stages of valve motion. An isogeometric biomechanical valve simulation is then performed on the fitted surfaces to show their direct applicability toward analysis. \textsc{ValveFit} enables automated patient-specific modeling with minimal manual intervention, paving the way for the future development of direct image-to-analysis platforms for clinical applications.
\end{abstract}


\end{frontmatter}
%%%%%%%%%%%%%%%%%%%%%%%%%%%%%%%%%%%%%%%%%%%%%%%
\section{Introduction}
\label{sec:intro}
Heart valve disease is one of the leading causes of heart failure, resulting in more than 25,000 deaths each year in the United States alone \cite{members2010heart, hinton2011heart, CDC}. This condition can manifest as \textit{stenosis}, which is the obstruction of blood flow; \textit{regurgitation}, improperly closing valves that result in backflow; or \textit{atresia}, where the valve is underdeveloped and lacks a natural opening for proper blood flow. If left untreated, heart valve disease can severely impair cardiac function; therefore, such medical conditions often require surgical repair or replacement of the valve. One of the severe cases of cardiac disease is congenital heart disease, in which children born with underdeveloped hearts often require several surgeries, further increasing the risk of surgical complications. In particular, patients with hypoplastic left heart syndrome (HLHS) undergo a surgical protocol that involves three successive palliative surgeries: the Norwood procedure, the Glenn procedure, and the Fontan procedure \cite{jayakumar2004cardiac, ohye2016current, gobergs2016hypoplastic, rai2019hypoplastic}. There is an increased risk of tricuspid valve regurgitation (TR) during these surgeries \cite{mccarthy2004tricuspid, Lee19fs, enriquez2019tricuspid, condello2021etiology, Ross23Tricu}, which, if not addressed in a timely manner, can cause right ventricle dilation or dysfunction.
 
In the past decade, advancements in non-invasive cardiac imaging methods have significantly improved the clinical diagnosis and treatment of heart valve diseases \cite{sun2014computational}. These advanced imaging modalities enable clinicians to generate time-resolved three-dimensional images non-invasively, allowing for preemptive and personalized treatment planning. In cardiac imaging, spatiotemporal echocardiography (4D echo) has emerged as an established non-invasive modality for the visualization of the heart valves, serving as an evaluation of the valve structure and function by providing anatomical and functional information \cite{lancellotti2013recommendations}. Despite these advances in non-invasive clinical imaging that enable accurate and earlier prediction of the onset of valve regurgitation and personalized treatment planning \cite{badano2009evaluation, muraru20193}, technical challenges persist in patient-specific computational modeling of the heart valves \cite{yan2021experimental}. The inherent noise, sparsity, and low resolution of imaging modalities such as echocardiography limit the accurate segmentation and visualization of tricuspid heart valves. Addressing these challenges requires the conversion of this low-fidelity imaging data to high-fidelity geometric representations that offer high-resolution surface visualization \cite{zhang2018geometric}. The solution also necessitates the integration of advanced surface fitting methods that can facilitate the automatic generation of patient-specific models, which can be directly used for computational biomechanical analysis \cite{sacks2017need}.

In this paper, we present a novel surface fitting framework that automatically generates patient-specific spline-based surface geometry of the tricuspid valve. The proposed method performs smooth and regular surface deformation of a template surface to fit the 3D point cloud obtained from echocardiographic image segmentation. We validate our framework using tricuspid valve surfaces obtained from simulation data, showing robustness toward different point cloud sampling and noise intensities. We then demonstrate the fitting results using point cloud data from real patients and the application of isogeometric analysis on the fitted surfaces ~\cite{johnson2021param}. The rest of the article is organized as follows: In Section~\ref{sec:literature}, we provide an overview of patient-specific computational modeling frameworks for heart valves. In Section~\ref{sec:imageacq}, we explain the acquisition of echocardiographic imaging data obtained from patients and the generation of patient-specific point clouds from echocardiographic data through image segmentation. The proposed \textsc{ValveFit} framework is introduced in Section~\ref{sec:methods}, and the surface fitting results along with a discussion on the performance are provided in Section~\ref{sec:results}. We also demonstrate that the surfaces generated by \textsc{ValveFit} are analysis-suitable by performing a simplified valve closure simulation. Finally, we provide concluding remarks and outline future directions in Section~\ref{sec:conclusion}.
%%%%%%%%%%%%%%%%%%%%%%%%%%%%%%%%%%%%%%%%%%%%%%%%

\section{Related Work on Patient-Specific Modeling of Heart Valves}\label{sec:literature}
The human heart consists of two atrioventricular valves, the tricuspid and mitral valves, and two semilunar valves, the pulmonary and aortic valves \cite{hinton2011heart}. The anatomical complexity of heart valves, which includes asymmetric geometry and thin and flexible leaflet structures, makes accurate patient-specific modeling particularly challenging \cite{sacks2009bioengineering, sacks2009biomechanics, toma2016high, gaidulis2022patient}. In addition, during the cardiac cycle, heart valves undergo drastic structural deformation, where the leaflets fold inward as they transition from the open to the closed configuration, further complicating the automatic surface fitting of the fast-moving valves from dynamic time-series images. 

Heart valve surface fitting methods play a crucial role in the functional and structural assessment of the valves by providing accurate geometric representations from medical imaging data \cite{haj2012general, miller2021implementation, toma2022clinical, bavo2017patient, bennati2023image, bennati2023turbulent, viola2022fsei, fedele2017patient, kaiser2021design, fumagalli2025reduced, govindarajan2024biomechanical, su2016cardiac}. Several frameworks have been presented in the realm of patient-specific modeling based on surface fitting of point clouds obtained from image segmentation \cite{Pouch14Fully, ginty2019dynamic}. For example, least squares optimization has been used to fit surfaces to point cloud segmentation data for the tricuspid valves \cite{mathur2022texas}. These methods rely on good-quality image segmentation data to generate analysis-suitable meshes. If the point cloud data from images is sparse, noisy, and incomplete, these existing methods will likely overfit the data, resulting in inaccurate patient-specific valve models that would require considerable manual intervention to fix the models in order to achieve good mesh quality for analysis \cite{votta2008mitral, stevanella2010finite}. Biomechanical models have been used to accurately fit the physiology of the heart valves \cite{mccarthy2004tricuspid, zhang2017towards, simonian2023patient, liu2023computational}. However, these methods do not impose stricter geometric regularization constraints to ensure good-quality surfaces during fast valve motion. On the other hand, machine learning-based approaches have been used to model the heart valve geometry directly from medical images, circumventing the image segmentation procedure \cite{munafo2024deep}. However, the accuracy of these methods is dependent on large data from diverse populations, which is typically not available.

B-spline and non-uniform rational B-spline (NURBS) surfaces provide smooth representations of the complex geometry of heart valves. In \cite{Xu18id, abbasi2020geometry, Pan24Param}, NURBS-based parameterization frameworks were developed to model aortic valves. Isogeometric analysis (IGA) \cite{Hughes05a, li2019isogeometric, wei2015truncated, wei2016extended, wei2017truncated, wei2017truncatedt, wei2018blended, casquero2020seamless}, which leverages spline-based geometric representations for direct analysis, has been successfully applied to patient-specific structural simulations of aortic valve closure, demonstrating high efficiency compared to traditional finite element (FE) methods \cite{Morganti15fy}. NURBS-based aortic valve models have also been employed in various structural analyses \cite{Zhang21Isoge, Zhang21bhv, Feng23Funct}, fluid-flow simulations \cite{Takizawa17Heart, Terahara20bw, Terahara22Compu}, and fluid--structure interaction studies \cite{Wu19fm, Johnson20ke, Xu21ci, Johnson22Effec}, illustrating their versatility and effectiveness in cardiovascular applications. In \cite{balu2019deep}, a deep learning-based framework was proposed for the design and analysis of surgical valves, highlighting the potential of integrating machine learning techniques with geometric modeling approaches. Additionally, in \cite{Sacks22Neura, Meyer25Neura}, a NURBS-based trileaflet heart valve model was integrated into a neural network finite element (NNFE) framework, enabling rapid simulation of valve closure with a significant speed-up compared to traditional FE solvers. In \cite{johnson2021param}, a generalized modeling and analysis framework for tricuspid valves was introduced, featuring a novel parameterization method for atrioventricular valve modeling. However, the parameterization-based frameworks fall short of capturing the intricate and complex features of the heart valves \textit{automatically} from images and are not generalizable to large patient cohorts. Although surface fitting methods of deforming a template B-spline or NURBS surface to subject-specific point clouds have been proposed \cite{aggarwal2013patient, Aggarwal16ed, prasad2022nurbs, moola2023thb, Gramling2024Vivo}, these approaches are not robust to complex changes in the valve surface across a large patient database or for mapping dynamic surfaces at different time points of valve motion during the cardiac cycle. Under large deformation, these methods can produce warped or intersected surfaces that would need manual changes to the mesh to make it analysis-suitable \cite{park2021heart}. 

Existing surface fitting approaches face several challenges that limit their clinical utility. Current approaches risk overfitting to noisy and sparse point cloud data and therefore rely on high-quality image segmentation, which is difficult to achieve with real echocardiographic images. Many surface fitting approaches lack proper geometric regularization constraints during deformation, especially during large valve motion, resulting in twisting of surface geometry, which is not analysis-suitable. These approaches require manual interventions to correct unsatisfactory meshes for analysis. To address these limitations, we propose the \textsc{ValveFit} framework, where we address the challenges of generalizable heart valve fitting by using a smooth B-spline surface fitting differentiable framework with novel regularization constraints. The \textsc{ValveFit} framework enables the automatic generation of complex valve geometry. In the present study, we also demonstrate the effectiveness of the framework for large-scale, patient-specific modeling, addressing the challenges with respect to inter-patient anatomical variability. Finally, we will also demonstrate the robustness of our method to generate smooth and analysis-suitable patient-specific valve surfaces considering various levels of Gaussian noise in point cloud data. The key features of \textsc{ValveFit} can be summarized as follows:
\begin{itemize}
    \item A B-spline surface fitting framework using automatic differentiation, enabling faster B-spline surface evaluations and gradient computation.
    \item Novel geometric regularization loss terms that ensure mesh regularity and prevent mesh distortion even during large deformation, while still accurately fitting the point clouds.
    \item Ensuring robustness toward noise, point cloud densities, and anatomical variability between patients.
    \item Automatic generation of patient-specific, analysis-suitable valve geometries without the need for manual interventions.
\end{itemize}
These innovations enable the \textsc{ValveFit} framework to reliably generate analysis-suitable valve geometries directly from noisy and sparse clinical data with minimal manual intervention, laying the foundation for the future development of direct image-to-analysis platforms for clinical applications.

\section{Image Acquisition and Segmentation}\label{sec:imageacq}
This study is performed using echocardiographic imaging data from the University of Oklahoma Health Sciences Center of retrospectively enrolled HLHS-afflicted patients before the Stage I (Norwood) palliation surgery (IRB approval \#14112). In brief, routine 4-dimensional (4D, i.e., $X,Y,Z$, time) full-volume echocardiographic data is collected for multiple cardiac cycles using a Philips EPIQ ultrasound machine (Philips, NV) equipped with a 5 MHz transthoracic matrix-array transducer. Patient echo imaging data are obtained through ECG-gated acquisitions over 6 to 8 consecutive heartbeats to ensure reproducible and reliable data \cite{Mah2021}. The echocardiogram images are converted into a Cartesian Digital Imaging and Communications in Medicine (3D DICOM) format in QLAB Cardiac Analysis software (Philips, NV).

\begin{figure}
\centering
\includegraphics[width=\textwidth]{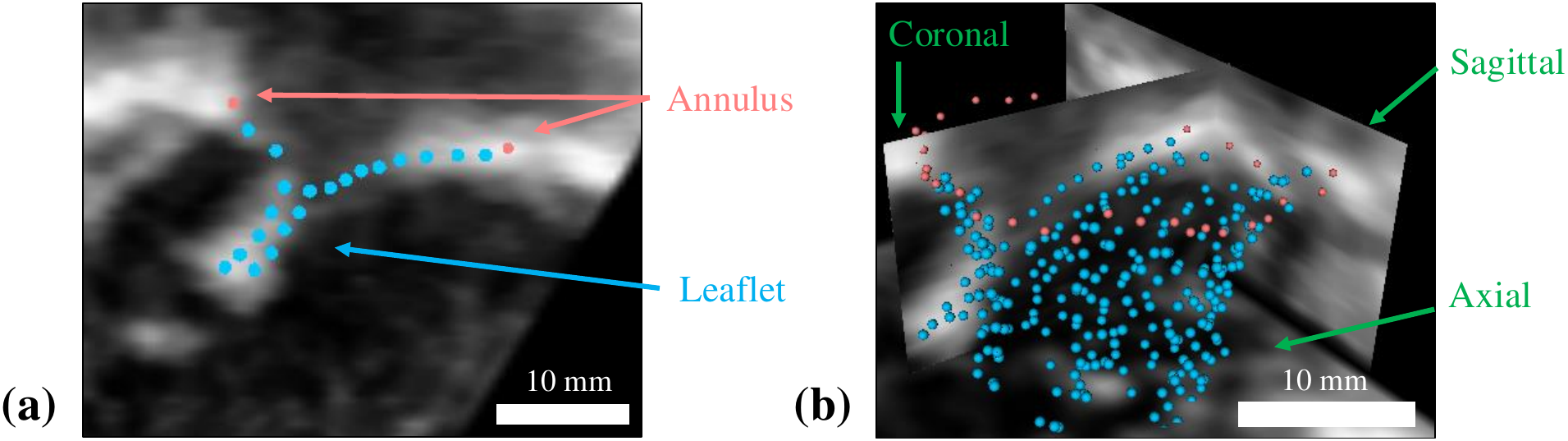}
\caption{Representative echocardiogram imaging data of a newborn with hypoplastic left heart syndrome (HLHS) and the segmented point cloud for the tricuspid valve (TV) apparatus: \textbf{(a)} coronal view and \textbf{(b)} 3D perspectives. Segmented points (red) for the TV annulus and segmented points (blue) for the three TV leaflets.}
\label{echo_image}
\end{figure}

Segmentations of the tricuspid valve are performed using 3D Slicer \cite{Slicer, fedorov20123d} at time instants between end-diastole (ED) and end-systole (ES) within the same cardiac cycle. First, the annulus is segmented by rotating the coronal image plane around the annulus in 10-degree increments and manually placing a point at the leaflet hinge point of each angular increment, which is identified by analyzing sequential time points in the echocardiogram sequence, to generate a point cloud consisting of 36 points (Fig.~\ref{echo_image}(a)). Then, points are added to manually trace the leaflets using the sagittal and coronal image frames with an approximate point density of 6.2 mm (Fig.~\ref{echo_image}(b)).

After the valves are segmented, the $X$, $Y$, and $Z$ coordinates of each leaflet are imported to \mbox{MATLAB} (MathWorks, USA), where an in-house program is used to classify each coordinate as belonging to the septal (SL), anterior (AL), or posterior leaflet (PL). Briefly, the $X$ and $Y$ coordinates of the commissure locations (i.e., the region between each leaflet), as manually identified using 3D Slicer, are first converted to polar coordinates $(r, \theta)$, with the septal-anterior commissure used as the reference angle (i.e., $\theta_{(SL-AL)}=0^{\degree}$). Using the values of $\theta$ at each commissure location, intervals are constructed for each leaflet for categorization of the leaflet segmentation data points. For example, a point $i$ was classified as belonging to the anterior leaflet if $\theta_{(SL-AL)} < \theta_{i}< \theta_{(AL-PL)}$.

%%%%%%%%%%%%%%%%%%%%%%%%%%%%%%%%%%%%%%%%%%%%%%%%%%%%%%%%%%%%%
\section{The \textsc{ValveFit} Framework}
\label{sec:methods}

\begin{figure}
    \centering
        \includegraphics[width=\linewidth]{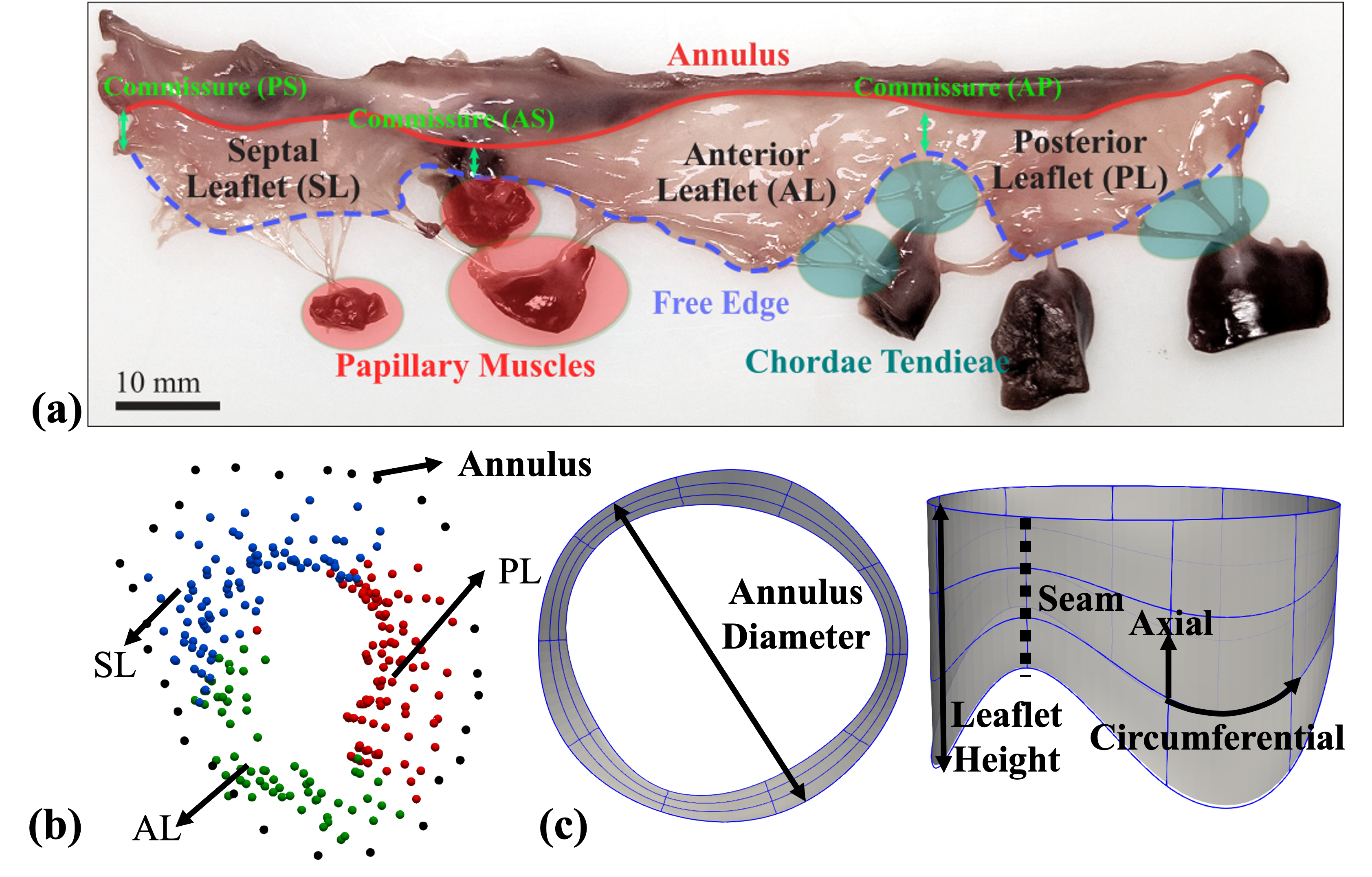}
        \caption{\textbf{(a)} Excised porcine tricuspid valve (TV) showing the three leaflets—septal (SL), anterior (AL), and posterior (PL)—along with their commissures. \textbf{(b)}  Point cloud representation of a patient’s TV annulus (black) and the three leaflets: SL (blue), AL (green), and PL (red). \textbf{(c)} Template geometry constructed using a single-patch periodic B-spline surface, shown in both top and side views. Idealized template geometry is generated with the annulus modeled as a circular ring and the leaflets with equal heights. Surface parameterization is carried out using knot vectors: $\mathcal{U} := \{0, 0, 0, 0.33, 0.66, 1, 1, 1\}$ along the axial direction and $\mathcal{V} := \{-0.3, -0.2, -0.1, 0, 0.1, 0.2, 0.3, 0.4, 0.5, 0.6, 0.7, 0.8, 0.9, 1, 1.1, 1.2, 1.3\}$ along the circumferential direction.}
        \label{fig:template}
\end{figure}

The proposed \textsc{ValveFit} framework generates an accurate and smooth surface representation of the tricuspid heart valve using B-splines. This is achieved by deforming a template B-spline surface to fit unstructured point clouds derived from the image segmentation of 4D echocardiographic data. The surface is deformed by moving control points based on the minimization of a loss function using gradient-descent optimization. The template geometry is constructed as a single patch B-spline surface, informed by the anatomical structure of the tricuspid valve, in order to ensure topological similarity. The anterior, posterior, and septal leaflets of the tricuspid heart valve (shown in Fig.~\ref{fig:template}(a)) are attached to a fibrous ring known as the annulus. These leaflets are anchored to the papillary muscles of the right ventricle through chordae tendineae, which prevent the leaflets from inverting during ventricular contraction \cite{hinton2011heart, miller2021implementation}.

\subsection{The Template Geometry}
\begin{figure}
    \centering
    \begin{subfigure}[b]{0.43\linewidth}
        \centering
        \includegraphics[height=.87\linewidth]{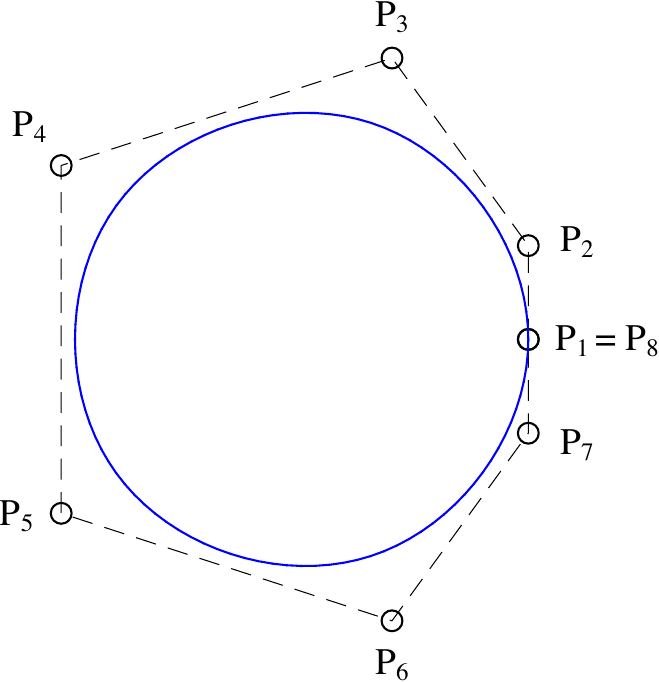}
        \caption{}
        \label{fig:unclamped_closed}
    \end{subfigure}
    \begin{subfigure}[b]{0.43\linewidth}
        \centering
        \includegraphics[height=.87\linewidth]{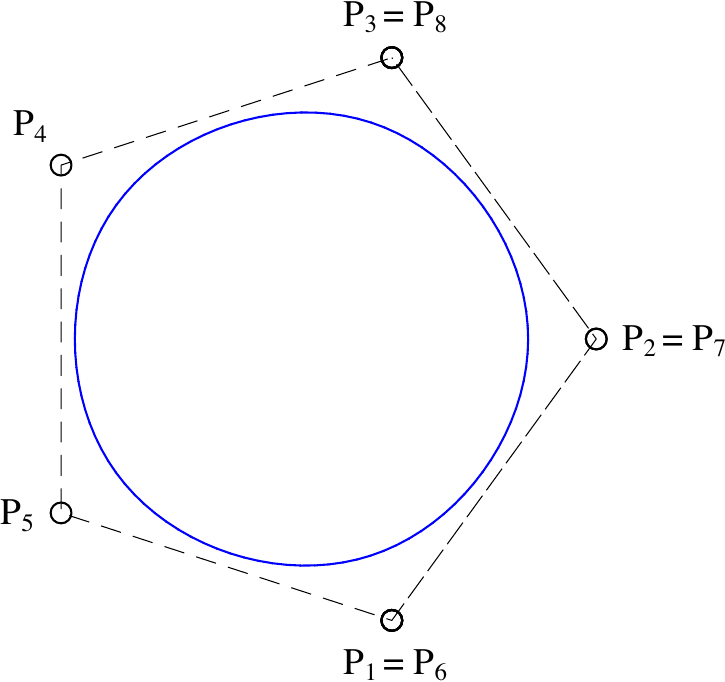}
        \caption{}
        \label{fig:clamped_closed}
    \end{subfigure}
    \caption{Comparison of \textbf{(a)} closed non-periodic and \textbf{(b)} closed periodic cubic B-spline curves. The periodic curve in \textbf{(b)} achieves $C^{2}$ continuity across the parametric boundary, whereas the non-periodic curve in \textbf{(a)} maintains only $C^{0}$ continuity.}
    \label{fig:combined_basis_curves}
\end{figure}

Due to its annular ring and leaflet configuration, the tricuspid valve topology is similar to that of a cylinder. In the template geometry, the annulus is modeled as a circular ring, and the leaflets are assigned equal heights to represent an idealized tricuspid valve in the open configuration during diastole (see Fig.~\ref{fig:template}(c)) \cite{johnson2021param, toma2022clinical}. Note that a cylindrical surface could be directly used as the template geometry due to its topological similarity to the tricuspid valve. However, we chose an idealized template TV surface to better capture the key anatomical features of the target point clouds shown in Fig.~\ref{fig:template}(b), including the three-leaflet topology and the characteristic annular shape.
\begin{figure}
    \centering
    \begin{subfigure}[]{0.49\textwidth}
        \centering
        \includegraphics[width=\linewidth]{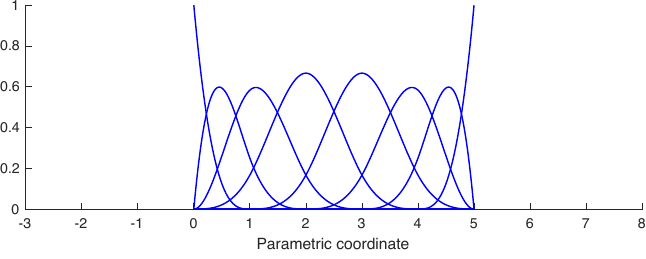}
        \caption{}
        \label{fig:unclamped_basis}
    \end{subfigure}    
    \hfill  
    \begin{subfigure}[]{0.49\textwidth}
        \centering
        \includegraphics[width=\linewidth]{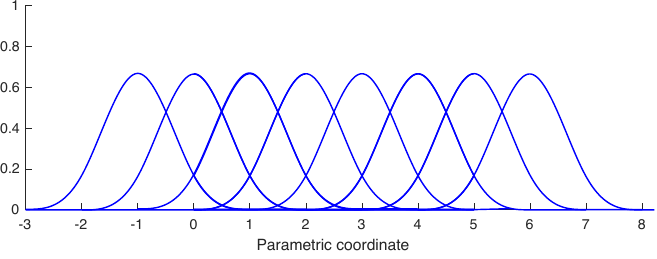}
        \caption{}
        \label{fig:clamped_basis}
    \end{subfigure}
    \caption{Comparison of cubic B-spline basis functions for \textbf{(a)} clamped and \textbf{(b)} unclamped knot vectors, defined as $\{0, 0, 0, 0, 1, 2, 3, 4, 5, 5, 5, 5\}$ and $\{-3, -2, -1, 0, 1, 2, 3, 4, 5, 6, 7, 8\}$, respectively.}
    \label{fig:combined_basis_fns}
\end{figure}
To enable higher continuity across the seam (Fig.~\ref{fig:template}(c)), the template geometry is constructed as a single patch B-spline surface using periodic B-splines. This B-spline surface undergoes non-linear deformation through loss function minimization, resulting in a patient-specific TV B-spline surface at the end of the fitting process. A B-spline surface ($\mathbf{S} \in \mathbb{R}^3$) can be described as 
\begin{equation}
    \mathbf{S}(u, v) = \sum_{i=1}^{n}\sum_{j=1}^{m} B_{i, p}(u)B_{j, q}(v)\mathbf{P}_{i, j} \, ,
    \label{eq:surface}
\end{equation} where $B_{i, p}(u)$ and $B_{j, q}(v)$ are univariate B-spline basis functions of degrees $p$ and $q$, defined on the knot vectors $\mathcal{U} := \{u_1, u_2, \dots, u_{n + p + 1}\}$ and $\mathcal{V} := \{v_1, v_2, \dots, v_{m+q+1}\}$, corresponding to the $u$- (circumferential) and $v$- (axial) parametric directions, respectively \cite{PiegTil97}. The surface point is given by $\mathbf{S}(u, v) = (x(u,v), \, y(u,v), \, z(u,v))$, where $x(u,v)$, $y(u,v)$, and $z(u,v)$ represent the Cartesian coordinates of the point.  Here, $n$ and $m$ denote the number of basis functions in each parametric direction, and $\mathbf{P}_{i,j} \in \mathbb{R}^3$ are the control points. The template geometry is constructed as a periodic B-spline surface with periodicity enforced along the circumferential direction. Accordingly, $B_{i, p}(u)$ is defined as an \emph{open non-periodic} B-spline using \emph{clamped} knot vectors, while $B_{j, q}(v)$ is a \emph{closed periodic} B-spline defined over \emph{unclamped} knot vectors. To create a closed B-spline curve, which can be either periodic or non-periodic, the control points at the two ends of the curve must be closed or overlapped. The distinction between a closed non-periodic and a closed periodic B-spline curve, along with their respective constructions, is shown in Figs.~\ref{fig:combined_basis_fns} and \ref{fig:combined_basis_curves}.  In Fig.~\ref{fig:combined_basis_fns}(a), cubic B-spline basis functions are defined over the \emph{clamped} knot vector: $\{0, 0, 0, 0, 1, 2, 3, 4, 5, 5, 5, 5\}$. As seen in Fig.~\ref{fig:combined_basis_curves}(a), overlapping the first and last control points results in a closed B-spline curve. However, the continuity at the parametric boundary is reduced to $C^{0}$ due to knot multiplicity at the ends of the knot vector. In Fig.~\ref{fig:combined_basis_fns}(b), cubic basis functions are defined over the \emph{unclamped} knot vector: $\{-3, -2, -1, 0, 1, 2, 3, 4, 5, 6, 7, 8\}$. By wrapping three control points at the two ends of the curve, we can achieve $C^{2}$ continuity in the interior as well as across the parametric boundary, as seen in Fig.~\ref{fig:combined_basis_curves}(b)~\cite{PiegTil97}.

\vspace{-0.5 em}

\subsection{The Loss Function}
\label{sec:loss-fn}
The proposed \textsc{ValveFit} framework performs iterative deformation through the minimization of a loss function using gradient descent optimization. To ensure that the result of our surface fitting framework is an analysis-suitable surface approximating the patient-specific TV geometry \cite{li2014analysis, wei2022analysis}, we construct a loss function which quantifies the shape dissimilarity between the template surface and the target point cloud along with surface mesh regularity. The loss function ($\mathcal{L}$), which consists of a shape fidelity term and a surface regularization term is defined as \begin{equation}
\mathcal{L}(\mathbf{S}, \mathbf{Q}) = \mathcal{E}_{\text{fid}}(\mathbf{S}, \mathbf{Q}) + \mathcal{E}_{\text{reg}}(\mathbf{S})\, ,
\label{eq:loss}
\end{equation} where $\mathcal{E}_{\text{fid}}(\mathbf{S}, \mathbf{Q})$ is the shape fidelity loss and $\mathcal{E}_{\text{reg}}(\mathbf{S})$ is the surface regularization energy. $\mathbf{S}$ denotes the template B-spline surface, and $\mathbf{Q}$ denotes the target point cloud.

\subsubsection{Shape Fidelity Loss} 
\label{sec:shape-fidelity}
The shape fidelity loss $\mathcal{E}_{\text{fid}}(\mathbf{S}, \mathbf{Q})$ ensures that the fitted surface closely approximates the target point cloud. It is defined as \begin{equation}
\mathcal{E}_{\text{fid}}(\mathbf{S}, \mathbf{Q}) = w_{\text{CD}} \, d_{\mathrm{CD}}(\mathbf{S}, \mathbf{Q}) + w_{\text{HD}} \, d_{\mathrm{HD}}(\mathbf{S}, \mathbf{Q}) + w_{\text{a}} \, d_\text{a}(\mathbf{A}_{\text{S}}, \mathbf{A}_{\text{Q}})\text{ ,}
\label{eq:data-fidelity}
\end{equation} Here, $d_{\mathrm{CD}}(\mathbf{S}, \mathbf{Q})$ is the one-sided Chamfer distance, $d_{\mathrm{HD}}(\mathbf{S}, \mathbf{Q})$ is the Hausdorff distance, and $d_{\mathrm{a}}(\mathbf{A}_{\text{S}}, \mathbf{A}_{\text{Q}})$ enforces the annulus constraint. The weights $w_\text{CD}$, $w_\text{HD}$, and $w_\text{a}$ control the influence of each term.

\begin{figure}[htb!]
    \centering
    \includegraphics[width=\linewidth]{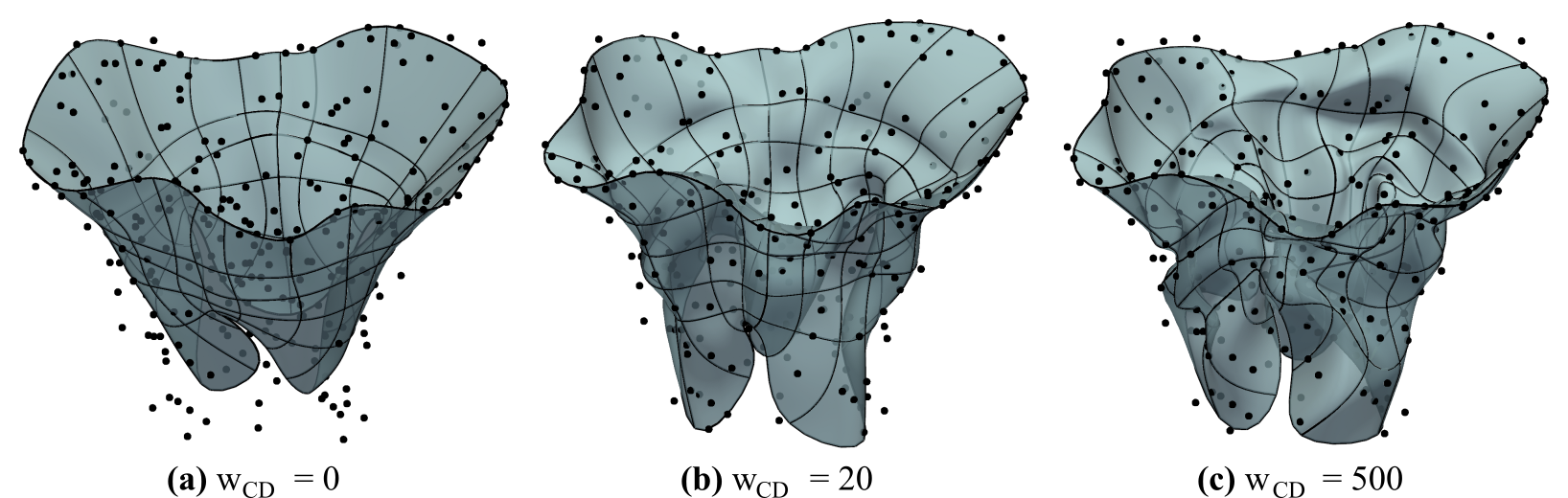}
    \caption{Effect of the Chamfer distance weight, $w_\text{CD}$, on surface fitting accuracy. \textbf{(a)} With no contribution from the Chamfer distance term ($w_\text{CD} = 0$), the surface fails to accurately fit to the target point cloud. \textbf{(b)} A moderate weight value ($w_\text{CD} = 20$) yields a well-fitted surface that accurately approximates the target point cloud. \textbf{(c)} An excessive large weight value ($w_\text{CD} = 500$) leads to overfitting, resulting in unnatural surface folding.}
    \label{fig:chamfer-distance}
\end{figure}

\textit{\textbf{Chamfer Distance}} ($d_\text{CD}$). It measures the average discrepancy between two point clouds by calculating the mean nearest-neighbor distance in both directions. Given a template point set $\mathbf{S}$ and a target point set $\mathbf{Q}$, the symmetric Chamfer distance, $d^{sym}_{\mathrm{CD}}(\mathbf{S}, \mathbf{Q})$, is defined as

\begin{equation}
 d^{sym}_{\mathrm{CD}}(\mathbf{S}, \mathbf{Q}) = \frac{1}{M}\sum_{k=1}^{M} \min_{l} \|\mathbf{s}_k - \mathbf{q}_l\|^2 + \frac{1}{N}\sum_{l=1}^{N} \min_{k} \|\mathbf{q}_l - \mathbf{s}_k\|^2.
 \label{eq:symmetric-chamfer}
\end{equation} Here $\mathbf{s}_k = \mathbf{S}(u_k,v_k)$, for $k = 1, ..., M$, denote points sampled from the template surface $\mathbf{S}$, and $\mathbf{q}_l\in \mathbf{Q}$, for $l = 1, ..., N$, denote points in the target point cloud $\mathbf{Q}$. The first term measures the mean nearest neighbor distance from the template surface points to the target point cloud ($\mathbf{S} \to \mathbf{Q}$), while the second term measures the mean nearest neighbor distance from the target point cloud to the template surface points ($\mathbf{Q} \to \mathbf{S}$). While the symmetric Chamfer distance is a comprehensive measure of bidirectional similarity, it can introduce artifacts when the target point cloud $\mathbf{Q}$ is sparse. The first term ($\mathbf{S} \to \mathbf{Q}$) forces every point on the dense template surface to align closely with one of the few sparse target points, which may cause the template surface to collapse or self-intersect in order to minimize this distance. To avoid these artifacts and ensure a more robust fit, we instead employ a one-sided Chamfer distance that considers only the second term, measuring the nearest neighbor distance from the target point cloud $\mathbf{Q}$ to the template surface $\mathbf{S}$. This is defined as

\begin{equation}
 d_{\mathrm{CD}}(\mathbf{S},\mathbf{Q}) = \frac{1}{N}\sum_{l=1}^{N} \min_{k} \|\mathbf{q}_l - \mathbf{s}_k\|^2.
 \label{eq:one-sided-chamfer-q-to-s}
\end{equation} Fig.~\ref{fig:chamfer-distance} demonstrates the effect of the Chamfer distance weight $w_\text{CD}$. A weight that is too low leads to underfitting, where the surface fails to accurately fit the target point cloud, while a weight value that is too high causes overfitting and introduces artifacts such as surface folding.

\begin{figure}[htb!]
    \centering
    \includegraphics[width=\linewidth]{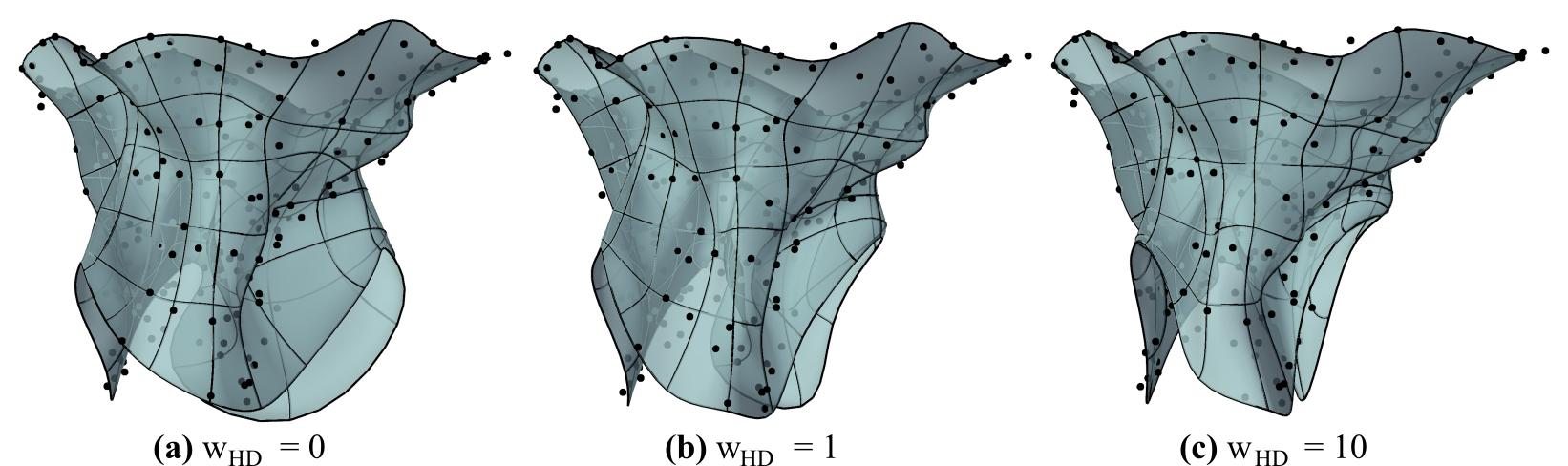}
    \caption{Effect of Hausdorff distance weight, $w_\text{HD}$, on constraining surface boundaries. \textbf{(a)} With $w_\text{HD} = 0$, the right-bottom boundary of the surface extends significantly beyond the target point cloud. \textbf{(b)} A moderate weight value ($w_\text{HD} = 1$) pulls the surface boundary closer, improving alignment with the target point cloud. \textbf{(c)} A higher weight value of $w_\text{HD} = 10$ enforces tight adherence of the surface boundary to the boundary target points.}
    \label{fig:hausdorff-distance}
\end{figure}

\textit{\textbf{Hausdorff Distance}} ($d_\text{HD}$). The Hausdorff distance term ensures that the maximum deviation between the template surface and the target point cloud is minimized, thus ensuring that even the most distant points, often located near the boundaries, are closely approximated. It is defined as

\begin{equation}
    d_{\mathrm{HD}}(\mathbf{S}, \mathbf{Q}) = \max \left( \max_k \, \min_l \|\mathbf{s}_k - \mathbf{q}_l\|, \max_l \, \min_k \|\mathbf{q}_l - \mathbf{s}_k\| \right).
    \label{eq:hausdorff}
\end{equation} The use of a one-sided Chamfer distance allows the template surface boundary to expand beyond the extent of the target point cloud without incurring a penalty, as it only enforces the proximity of the target points to the template surface. In contrast, the Hausdorff distance term explicitly penalizes the maximum deviation between the surface and the point cloud, thus constraining surface growth beyond the local neighborhood of the target point cloud. This effect is illustrated in Fig.~\ref{fig:hausdorff-distance} by varying the weight ($w_\text{HD}$) associated with the Hausdorff distance in Eq.~\eqref{eq:hausdorff}, where an enlarged template surface is fitted to the target point cloud.

\begin{figure}[htb!]
    \centering
    \includegraphics[width=\linewidth]{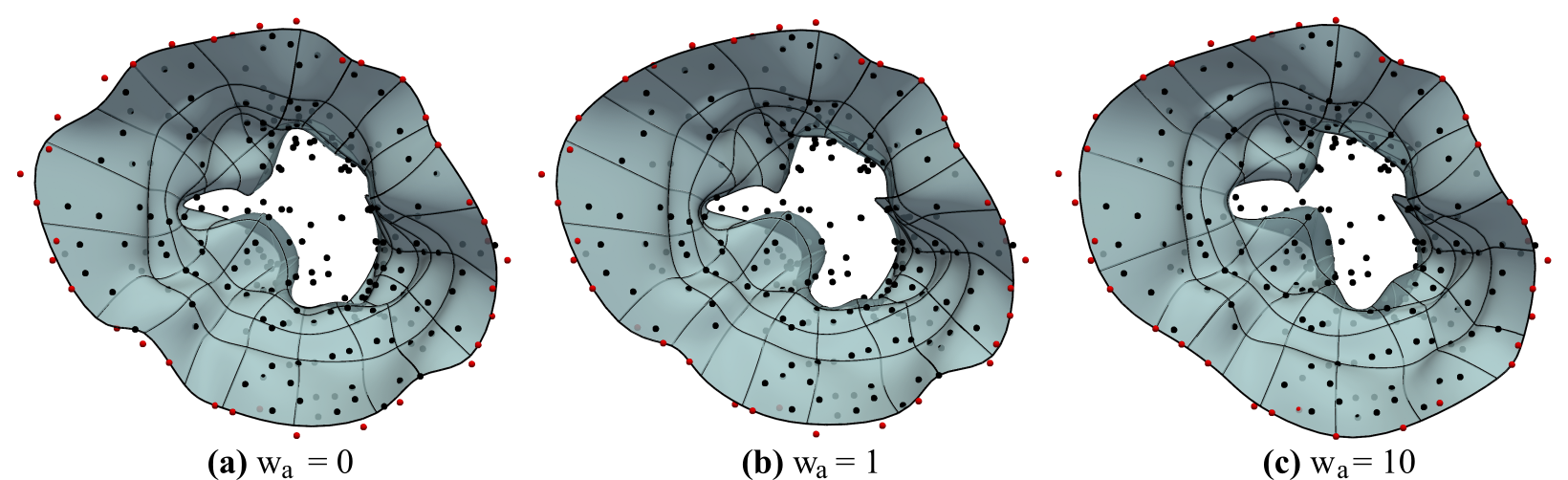}
    \caption{Effect of the weight, $w_a$, for fitting to the target annulus points (red). \textbf{(a)} Weight of zero; \textbf{(b)} a moderate weight aligns the surface boundary with the target annulus points; and \textbf{(c)} a higher weight enforces a tighter conformation.}
    \label{fig:annulus}
\end{figure}

\textit{\textbf{Annulus Constraint}} ($d_\text{a}$). The term $d_{\mathrm{a}}(\mathbf{A}_{\text{S}}, \mathbf{A}_{\text{Q}})$ calculates the symmetric Chamfer distance (Eq.~\eqref{eq:symmetric-chamfer}) between the parametric boundary on the template surface closest to the annulus $\mathbf{A}_{\text{S}} := \mathbf{S}(u=u_0)$  and the segmented points  corresponding to the annulus in the target point cloud $\mathbf{A}_{\text{Q}}$. The effect of this constraint is demonstrated in Fig.~\ref{fig:annulus} by varying its corresponding weight, $w_a$. When the weight is zero (Fig.~\ref{fig:annulus}\textbf{(a)}), the annulus constraint is disabled. This results in a loose conformation that does not explicitly capture the target annulus points. As $w_a$ increases (Fig.~\ref{fig:annulus}\textbf{(b)-(c)}), this term enforces a better alignment, pulling the template surface boundary to achieve a tighter fit with the target annulus points.

\subsubsection{Surface Regularization Energy}
\label{sec:regularization}
The surface regularization energy $\mathcal{E}_{\text{reg}}(\mathbf{S})$ imposes smoothness and regularity constraints on the deforming surface to obtain a smooth parameterization in the physical domain and prevent element distortions or self-intersections. It  comprises three terms: \begin{equation}
\mathcal{E}_{\text{reg}}(\mathbf{S}) = w_{\text{orth}} \, R_{\text{orth}}(\mathbf{S}) + w_\text{TPE} \, R_\text{TPE}(\mathbf{S}) + w_{\text{norm}} \, R_{\text{norm}}(\mathbf{S}).
\label{eq:regularization}
\end{equation} The weights $w_{\text{orth}}$, $w_{\text{norm}}$, and $w_\text{TPE}$ allow us to control the influence of each of the above three terms in Eq.~\eqref{eq:regularization}, allowing for adjustments based on specific mesh quality requirements.

\begin{figure}[htb!]
    \centering
    \includegraphics[width=1.0\linewidth]{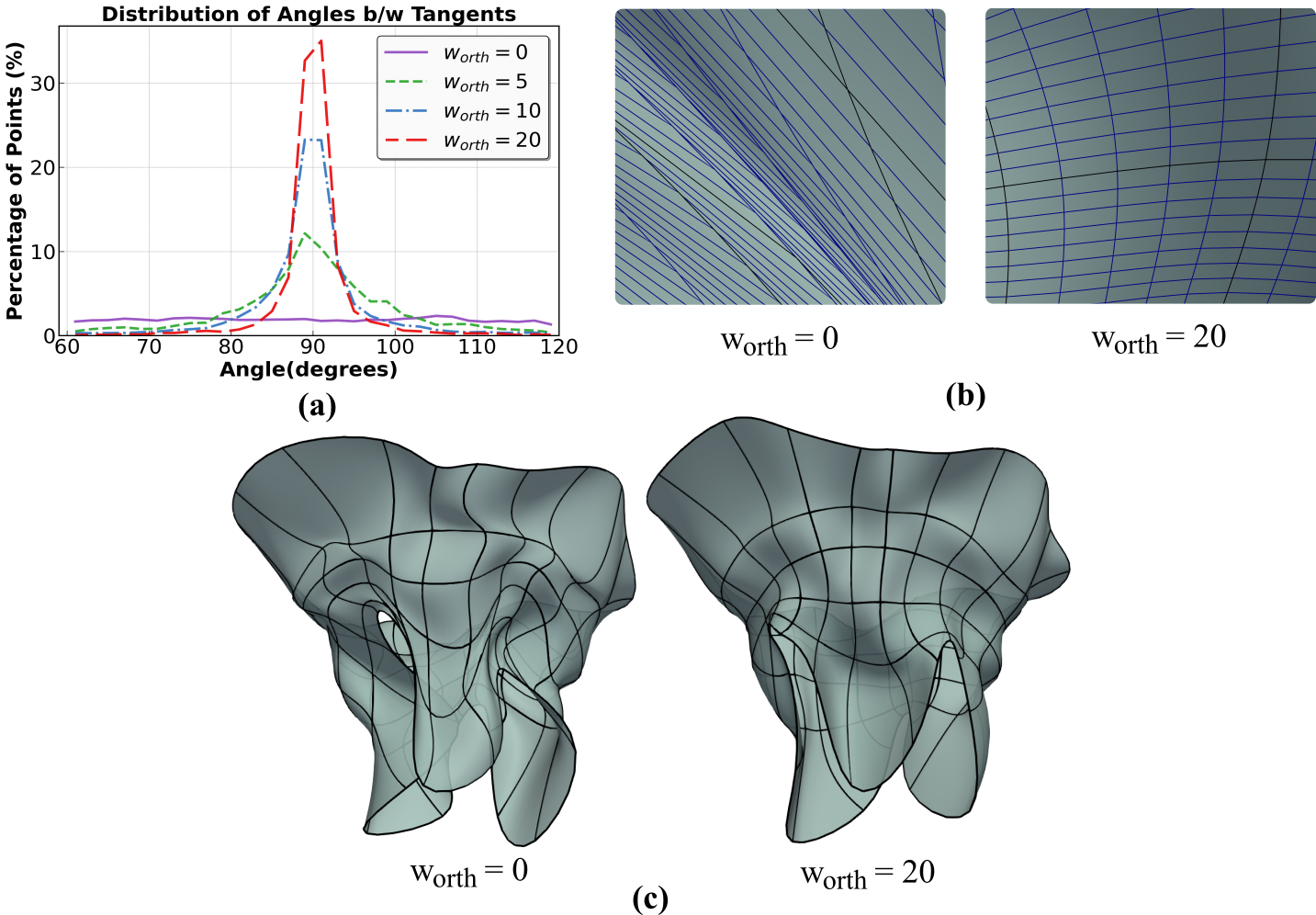}
    \caption{The effect of varying the orthogonality weight, $w_\text{orth}$. \textbf{(a)}~Distribution of the angle between tangent vectors, showing improved orthogonality with higher weight values. A visualization of the \textbf{(b)}~surface elements and \textbf{(c)}~fitted surface at $w_\text{orth}=0$ and  $w_\text{orth}=20$.}

    \label{fig:w_orth_study}
\end{figure}

\textit{\textbf{Tangent Orthogonality Energy}} ($R_\text{orth}$). To enforce bijectivity during surface deformation, we need to avoid self-intersections and preserve surface mesh element quality throughout the transformation \cite{gain2002preventing}. In particular, we ensure good surface mesh quality by penalizing extreme mesh distortion in each element. This is achieved by introducing an orthogonality regularization term, $R_{\text{orth}}(\mathbf{S})$, which ensures that the two parametric lines in the physical domain are locally orthogonal, thus maintaining surface regularity during deformation. This term is defined as
\begin{equation}
R_{\text{orth}}(\mathbf{S}) = \frac{1}{M} \sum_{k=1}^{M} \left| \frac{\mathbf{t}_k^u \cdot \mathbf{t}_k^v}{\left\| \mathbf{t}_k^u \right\| \left\| \mathbf{t}_k^v \right\|} \right|,
\label{eq:orthogonality-term}
\end{equation}
where $\mathbf{t}_k^u=\mathbf{t}^u(u_k, v_k)$ and $\mathbf{t}_k^v=\mathbf{t}^v(u_k, v_k)$ denote the tangent vectors evaluated at each sampled parametric point $(u_k, v_k)$, and $M$ is the number of sampled points in the parametric domain. Minimizing this energy results in an orthogonal parameterization that is free of local shearing. This angle preserving quality is critical for ensuring an analysis-suitable surface. By calculating the cosine of the angle between the local tangent vectors at each parametric point, $R_{\text{orth}}$ term penalizes deviations from right angles, thus avoiding excessive skewing or twisting of surface elements during deformation. The influence of this term is controlled by a weight parameter, $w_\text{orth}$, which introduces a trade-off between geometric regularity and data-fitting accuracy. As illustrated in Fig.~\ref{fig:w_orth_study}, setting $w_\text{orth}$ to zero allows the surface to deform freely but can result in skewed mesh elements. Conversely, a higher weight value produces a smoother, more regular surface with a higher percentage of surface points with near locally orthogonal tangent vectors. $R_\text{orth}$ also ensures that the area of each mesh element remains positive during deformation, thus, also ensuring that the surface has no local degeneracies (it does not collapse to a point or line).

\begin{figure}[htb!]
    \centering
    \includegraphics[width=\linewidth]{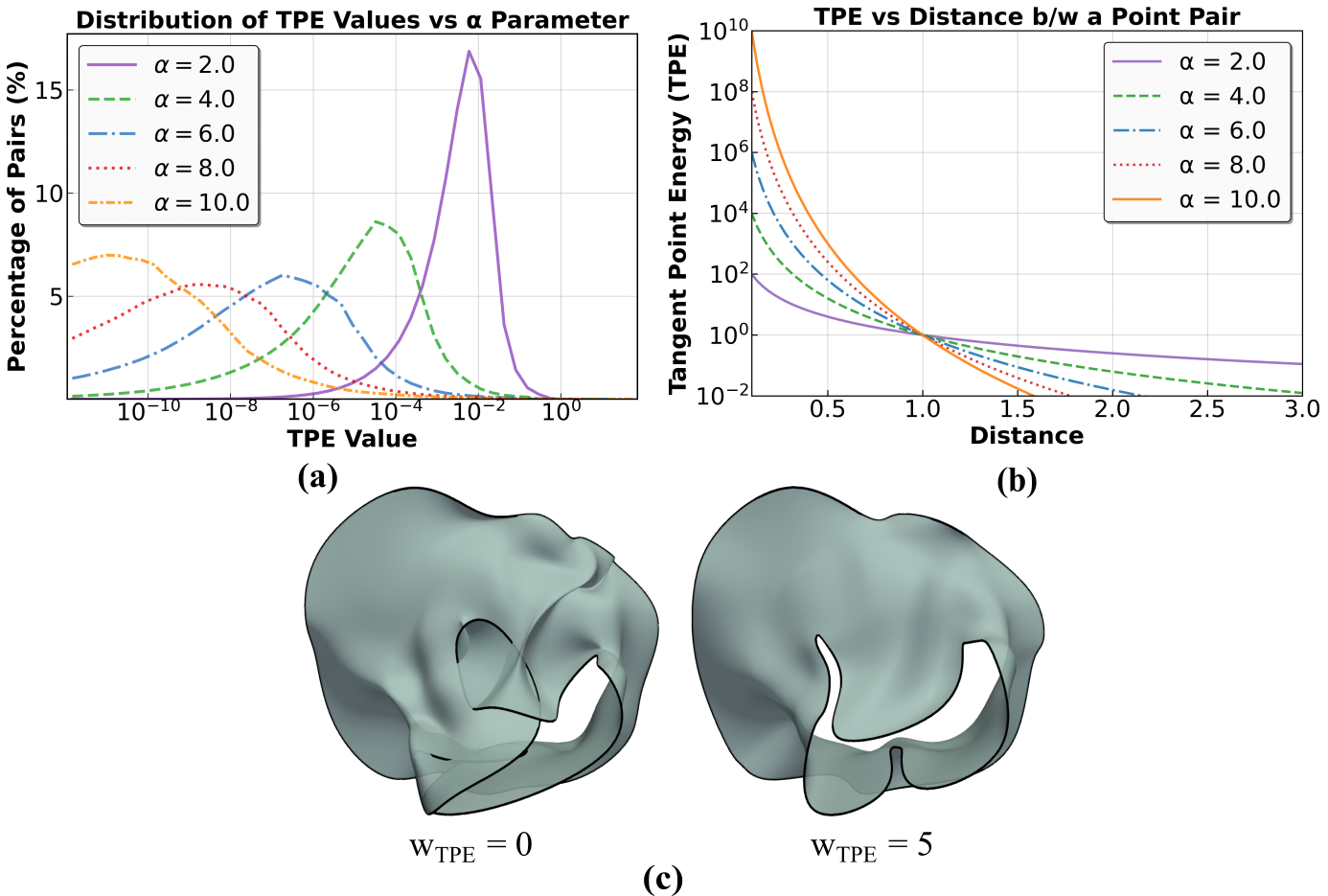}
    \caption{Effect of varying the Tangent Point Energy (TPE) parameters, $\alpha$ and $w_\text{TPE}$. \textbf{(a)}~Distribution of TPE values across all surface point pairs for different $\alpha$ values. \textbf{(b)}~TPE versus inter-point distance for different $\alpha$ values. \textbf{(c)}~Fitted surface with $w_\text{TPE}=0$, resulting in self-intersection, and with $w_\text{TPE}=5.0$, where TPE effectively prevents collisions.}
    \label{fig:TPE}
\end{figure}

\textbf{\textit{Tangent Point Energy}} ($R_\text{TPE}$). While the tangent orthogonality energy preserves local element shape quality, it does not prevent global self-intersections, such as collisions between leaflets during coaptation. To address this, we incorporate a tangent-point energy (TPE) term, $R_\text{TPE}(\mathbf{S})$, which introduces a repulsive potential that imposes an infinite energy barrier between any two surface points that approach each other \cite{Strzelecki2013, SaScRu24}. This term penalizes the most closest interaction between all  points on the surface, a strategy shown empirically to be more effective at preventing self-collisions than averaging over all interactions. This term is defined as
\begin{equation}
R_\text{TPE}(\mathbf{S}) = \max_{k \neq l} \left( \frac{\left| \mathbf{n}_k \cdot (\mathbf{s}_k - \mathbf{s}_l)  \right|^\alpha}{\left\| \mathbf{s}_k - \mathbf{s}_l \right\|^{2\alpha}} \right),
\label{eq:tpe-term}
\end{equation}
where $\mathbf{s}_k = \mathbf{S}(u_k, v_k)$ and $\mathbf{s}_l = \mathbf{S}(u_l, v_l)$ are any two distinct points on the surface, $ \mathbf{n}_k$ is the unit normal at $\mathbf{s}_k$, and $\alpha$ controls the strength and decay rate of the repulsive interaction force.

Fig.~\ref{fig:TPE} illustrates the influence of $\alpha$ and the weight $w_\text{TPE}$ on the TPE. $\alpha$ modulates the interaction range and sharpness of the repulsive potential. As shown in Fig.~\ref{fig:TPE}(b), higher values of $\alpha$ produce a steep, short-range repulsion, while lower values result in a gentler, long-range force. This behavior directly affects which point pairs are penalized; smaller $\alpha$ values distribute the repulsive energy across a broader set of pairs, whereas larger $\alpha$ values concentrate the penalty on a small subset of the closest set of pairs, as seen in Fig.~\ref{fig:TPE}(a). The accuracy of the TPE evaluation depends on the surface sampling strategy. Instead of carrying out adaptive quadrature \cite{SaScRu24}, we adopt a uniform sampling approach for computational efficiency. This method has proven sufficient to prevent self-intersections for the valve geometries considered in this study. The parameter value $\alpha=4.0$ is selected as it provides an effective balance between the repulsive force strength and its interaction range. As seen in Fig.~\ref{fig:TPE}(a) $\alpha=4.0$ focuses the penalty on a moderate group of high-risk point pairs more effectively than $\alpha=2.0$, while avoiding the narrow and unstable energy landscape associated with higher $\alpha$ values. The resulting $\alpha=4.0$ curve in Fig.~\ref{fig:TPE}(b) exhibits a smoother and more progressive onset of the repulsive contact force, resulting in stable convergence during optimization. We thus use $\alpha=4.0$ for all the examples in Section \ref{sec:results}. In Fig.~\ref{fig:TPE}(c), we showcase how by setting $w_\text{TPE}=5.0$, we can avoid leaflet collisions during coaptation by introducing the tangent point energy term.

\begin{figure}
    \centering
    \includegraphics[width=\linewidth]{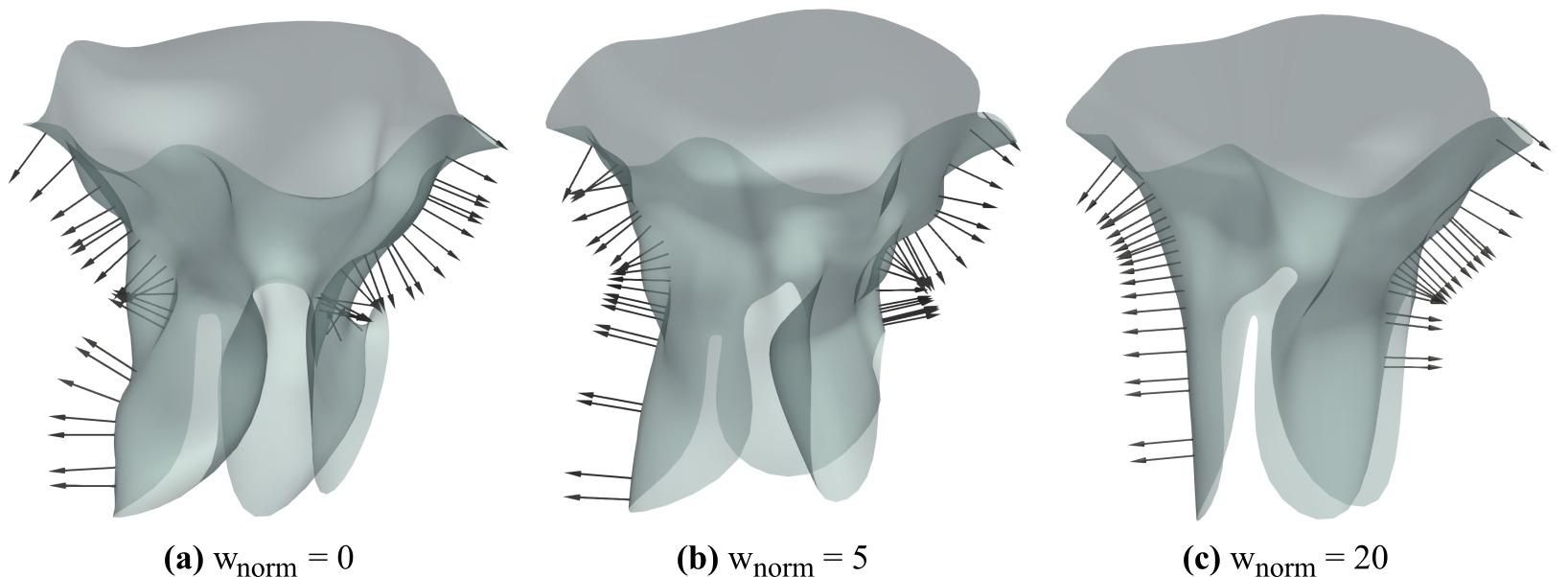}
    \caption{Effect of the normal deviation weight, $w_\text{norm}$, on surface smoothness and normal vector orientation. \textbf{(a)}~Without regularization ($w_\text{norm}=0$), the normal vectors intersect each other, particularly in the bottom right corner, resulting in the folding of the surface. \textbf{(b)}~A moderate weight ($w_\text{norm}=5.0$) reduces surface folding by enhancing consistent orientation of the normals. \textbf{(c)}~A high weight ($w_\text{norm}=20.0$) produces a regular surface with a well-aligned field of normal vectors.}
    \label{fig:w-norm}
\end{figure}

\textbf{\textit{Normal Deviation Energy}} ($R_\text{norm}$). The surface fitting process can be sensitive to noise in the target point cloud, often resulting in undesirable folding and self-intersections. To prevent this, we introduce a normal deviation energy term that reduces abrupt changes in the surface orientation. Specifically, by controlling the variation of the vertical ($z$) component of the unit normal vector $\mathbf{n}$, this term penalizes excessive folding, bending and self-intersections of the template surface. The normal deviation energy is defined as 
\begin{equation}
R_{\text{norm}}(\mathbf{S}) = \max_{1 \leq k \leq M} \left| \mathbf{n}_k^z - \bar{\mathbf{n}}^z \right|,
\label{eq:normal-variation}
\end{equation}
where $\mathbf{n}_k^z$ denotes the $z$-component of the unit normal vector $\mathbf{n}$ evaluated at the parametric coordinate $(u_k, v_k)$, and $\bar{\mathbf{n}}^z$ is the mean value of the $z$-component of the normal vectors across all the surface points. By minimizing this term, abrupt changes in surface orientation are prevented, suppressing excessive folding of the surface that may arise from fitting noisy point clouds. As demonstrated in Fig.~\ref{fig:w-norm}, increasing the weight $w_\text{norm}$ visibly regularizes the surface by enforcing a consistent alignment of the normal vectors and thus a more smoother geometry. Additionally, this term can be used as a soft constraint to control the shape of the surface, particularly in regions with sparse point cloud coverage, where data-driven guidance for fitting is difficult.

We can efficiently evaluate the tangent and normal vectors required for computing the terms in Eq.~\eqref{eq:regularization} by using the automatic differentiation capabilities of JAX \cite{jax2018github}. We evaluate the partial derivatives of the surface coordinates ($x$, $y$, and $z$) with respect to the parametric coordinates $u$ and $v$ by evaluating the Jacobian matrix ($J$) of the B-spline surface parameterization as
\begin{equation}
\renewcommand{\arraystretch}{1.5}
J = \begin{bmatrix}
\dfrac{\partial x}{\partial u} & \dfrac{\partial y}{\partial u} & \dfrac{\partial z}{\partial u} \\
\dfrac{\partial x}{\partial v} & \dfrac{\partial y}{\partial v} & \dfrac{\partial z}{\partial v}
\end{bmatrix}^{T}.
\label{eq:jacobian}
\end{equation} 
The tangent vectors ($\mathbf{t}^{u}$ and $\mathbf{t}^{v}$) and the unit surface normal vector $\mathbf{n}$ are derived from the Jacobian matrix as follows:
\begin{align}
\mathbf{t}^u(u, v) &= \frac{\partial \mathbf{S}}{\partial u} = \begin{bmatrix}
\dfrac{\partial x}{\partial u} \
\dfrac{\partial y}{\partial u} \
\dfrac{\partial z}{\partial u}
\end{bmatrix}^T \label{eq:tangent-u} \,, \\
\mathbf{t}^v(u, v) &= \frac{\partial \mathbf{S}}{\partial v} = \begin{bmatrix}
\dfrac{\partial x}{\partial v} \
\dfrac{\partial y}{\partial v} \
\dfrac{\partial z}{\partial v}
\end{bmatrix}^T \label{eq:tangent-v} \, , \\
\mathbf{n}(u, v) &= \frac{\mathbf{t}^u(u, v) \times \mathbf{t}^v(u, v)}{\left\| \mathbf{t}^u(u, v) \times \mathbf{t}^v(u, v) \right\|}.
\label{eq:normal-vector}
\end{align}

\subsection{Surface Fitting}
\label{sec:opt}
To accurately fit the template B-spline surface to target point clouds obtained from image segmentation, we use differentiable programming during surface fitting \cite{prasad2022nurbs, moola2023thb}. In the template B-spline surface (Eq.~\eqref{eq:surface}), the parametric coordinates $u$ and $v$ are fixed during optimization, and the control point positions $\mathbf{P}$ are optimized to minimize the loss function. The gradient of the loss function with respect to $\mathbf{P}$ is given as  
\begin{equation}
    \nabla_\mathbf{P}\mathcal{L} = \frac{\partial \mathcal{L}}{\partial \mathbf{S}} \, \frac{\partial \mathbf{S}}{\partial \mathbf{P}} \text{.}
    \label{eq:chain-rule}
\end{equation} 
The partial derivative of the surface with respect to each control point $\mathbf{P}_{i,j}$ is given by
\begin{equation}
    \frac{\partial \mathbf{S}}{\partial \mathbf{P}_{i,j}} = B_{i,p}(u) B_{j,q}(v) \text{.}\label{eq:bspline_derivative}
\end{equation} 
 During loss minimization, calculating the gradients of the loss function (Eq.~\eqref{eq:loss}) with respect to the surface $\mathbf{S}$ is more complex than computing the derivatives in Eq.~\eqref{eq:bspline_derivative}.  For example, the gradient of the Chamfer distance term is given as
\begin{equation}
    \frac{\partial d_\text{CD}}{\partial \mathbf{S}} = \frac{\partial}{\partial \mathbf{S}} \left(\frac{1}{N}\sum_{l=1}^N \min_{k} \|\mathbf{s}_k - \mathbf{q}_l\|\right)\text{,}
\label{eq:chamfer_gradient}
\end{equation} where $\mathbf{s}_k$ are points sampled from the B-spline surface and $\mathbf{q}_l$ are points from the target point cloud. Manual gradient evaluation is laborious and ineffective since it would be necessary to derive and implement the gradients for each term in the loss function, which would involve significant bookkeeping and be prone to error. This complexity emphasizes how crucial it is to use automatic differentiation in our surface fitting framework in order to compute gradients effectively and precisely. The gradient of the loss function with respect to all control points is given by
\begin{align}
\quad & \nabla_\mathbf{P} \mathcal{L} = \left[\frac{\partial \mathcal{L}}{\partial \mathbf{P}_{1,1}}, \frac{\partial \mathcal{L}}{\partial \mathbf{P}_{1,2}}, \ldots, \frac{\partial \mathcal{L}}{\partial \mathbf{P}_{n,m}}\right]^T,
\end{align} and each control point is updated iteratively as
\begin{align}
    \quad & \mathbf{P}_{i,j}^{(t+1)} = \mathbf{P}_{i,j}^{(t)} - \delta \frac{\partial \mathcal{L}}{\partial \mathbf{P}_{i,j}^{(t)}},
    \label{eq:control-point-update}
\end{align} where $\delta$ is the step-size and $t$ denotes the time step. We have utilized a step-size of $0.001$, determined empirically based on the maximum number of iterations and rate of convergence. The optimization of control points $\mathbf{P}$ to minimize the loss function is performed using the Adam optimizer \cite{kingma2017adam}.

The surface fitting process begins by initializing a template B-spline surface ($\mathbf{S}_0$) with predefined basis functions $(\mathbf{B}_u, \mathbf{B}_v)$ and initial control points ($\mathbf{P}_0$) at the start of the optimization process. The target data consist of a sequence of point clouds $\{\mathbf{Q}_1, \mathbf{Q}_2, \dots , \mathbf{Q}_n\}$, each segmented from echocardiographic images at various stages of cardiac cycle. These point clouds are ordered from the fully open state to the fully closed state. 

For each point cloud in the sequence, the fitting process starts by initializing the control points of the surface with those from the previous stage. For the first point cloud $\mathbf{Q}_1$, the process begins with the control points of initial template $\mathbf{P}_0$. To provide a good initial alignment, an affine transformation is applied to the surface so that it matches the global position, scale, and orientation of the target point cloud. The template surface is iteratively deformed to match the target point cloud driven by the gradient of the loss function (Eq.~\eqref{eq:control-point-update}). A step size $\delta=0.001$ and a maximum of $t_{max}=10,000$ iterations were employed to ensure stable convergence and accurate surface fitting.

Upon convergence for a given point cloud, the optimized control points are stored. These control points are then used to initialize the fitting process for the subsequent point cloud. Repeating this procedure across the entire sequence yields a series of control points $\{\mathbf{P}_1, \mathbf{P}_2, \dots, \mathbf{P}_n\}$ corresponding to the fitted surfaces $\{\mathbf{S}_1, \mathbf{S}_2, \dots, \mathbf{S}_n\}$. The program outlined in Algorithm~\ref{alg:surface_fitting} is implemented utilizing JAX arrays, which runs array operations on both GPU and CPU.

\begin{algorithm}[!htb]
\caption{\textsc{ValveFit}: Patient-Specific Surface Fitting for Tricuspid Valve Geometry}
\label{alg:surface_fitting}

\textbf{Input:} Template geometry with predefined parametrization: basis functions $(\mathbf{B}_u, \mathbf{B}_v)$ and control points $\mathbf{P}_0$; sequence of target point clouds $\{\mathbf{Q}_1, \mathbf{Q}_2, \dots, \mathbf{Q}_n\}$; maximum number of iterations $t_\text{max}$; step-size $\delta$; weights $\{w_\text{CD}, w_\text{HD}, w_\text{a}, w_\text{orth}, w_\text{norm}\}$

\textbf{Output:} Control points of fitted surfaces $\{\mathbf{P}_1, \mathbf{P}_2, \dots, \mathbf{P}_n\}$

\begin{algorithmic}[1]
    \FOR{$r = 1$ to $n$}
        \STATE Set: $\mathbf{P}_r \gets \mathbf{P}_{r-1}$ \hfill // Use previous fit as new template
        \STATE Apply affine transformation to align surface $\mathbf{S}_r$ to point cloud $\mathbf{Q}_r$: update $\mathbf{P}_r$
        \FOR{$t = 1$ to $t_\text{max}$}
            \STATE Compute surface points: $\mathbf{S}_r^{(t)} \gets \mathbf{B}_u \otimes \mathbf{B}_v \cdot \mathbf{P}_r^{(t)}$
            \STATE Evaluate loss: $\mathcal{L}(\mathbf{S}_r^{(t)}, \mathbf{Q}_r)$ with weights $\{w_\text{CD}, w_\text{HD}, w_\text{a}, w_\text{orth}, w_\text{norm}\}$ \hfill // Eq.~\eqref{eq:loss}
            \STATE Compute the gradient of the loss function with respect to control points: $\nabla \mathbf{P}_r^{(t)} \gets \frac{\partial \mathcal{L}}{\partial \mathbf{P}_r^{(t)}}$
            \STATE Update control points: $\mathbf{P}_r^{(t+1)} \gets \mathbf{P}_r^{(t)} - \delta \cdot \nabla \mathbf{P}_r^{(t)}$
        \ENDFOR
        \STATE Save final control points: $\mathbf{P}_r \gets \mathbf{P}_r^{(t_\text{max})}$
    \ENDFOR
\end{algorithmic}
\end{algorithm}
%%%%%%%%%%%%%%%%%%%%%%%%%%%%%%%%%%%%%%%%%%%%%%%%%%%%%%%%%%%%%%%%%%%%%%%%%%%%%

\section{Results}
\label{sec:results}
In this section, we demonstrate the performance of the \textsc{ValveFit} framework in accurately fitting B-spline surfaces to point cloud data. A validation study employing synthetic point clouds generated from structural simulations of heart valve motion is presented in Section~\ref{subsec:validation-study}. Since these point clouds are uniform and noise-free, they are used as the ground truth for TV geometry. Using these point clouds, we further investigate the effect of point cloud density (Section~\ref{subsec:density-study}) and noise (Section~\ref{subsec:noise-study}) on fitting accuracy. Finally, in Section~\ref{subsec:patient-data}, we apply the framework towards echocardiographic data acquired at multiple stages in the cardiac cycle to demonstrate its utility in generating patient-specific valve geometries.

To quantitatively evaluate the surface fitting accuracy, we introduce a dimensionless error metric, the scaled Nearest Neighbor Distance (sNND). While the Nearest Neighbor Distance (NND) measures the minimum distance from each point in the point cloud to the fitted surface, it is scale-dependent, varying with the size of the geometry. To address this limitation, sNND normalizes this nearest neighbor distance by the square root of the surface area, allowing for consistent comparison across valve geometries of different anatomical sizes in patients. For each point $\mathbf{q}_{l}$ in the target point cloud, where $l = 1, \cdots, N$, the sNND to the fitted surface of area $\mathcal{A}(\mathbf{S})$ is defined as 
\begin{equation}
    \text{sNND}(\mathbf{q}_l, \mathbf{S}) = \frac{\underset{k}{\min} \|\mathbf{q}_{l} - \mathbf{s}_k\|}{\sqrt{\mathcal{A}(\mathbf{S})}}\text{ .}
    \label{eq:snnd}
\end{equation}

The following results were generated on a hardware configuration consisting of an $\text{Intel}^\text{\textregistered}$ $\text{Core}^\text{TM}$ i9-14900KF CPU, 64 GB RAM, and an $\text{Nvidia}^\text{\textregistered}$ GeForce $\text{RTX}^\text{TM}$ 4090 GPU.

\subsection{Validation Study with Simulation Data}
\label{subsec:validation-study}
Point clouds obtained from image segmentation of echocardiogram images are often sparse, noisy, partial, and may not accurately represent the tricuspid valve surface. Ascertaining the accuracy of the fitted geometry using error metrics such as sNND, with respect to noisy point clouds, is not a reliable measure of the robustness of the fit, as the target point cloud itself may not accurately represent the underlying geometry. A very low sNND error for extremely noisy point clouds may in fact indicate that the fitted surface is overfitting the noise, resulting in an inaccurate representation of the underlying geometry. The resulting geometry is generally less smooth, of poor quality, and thus requires considerable manual effort to make it analysis-suitable. Therefore, we first validate our surface fitting algorithm by fitting B-spline surfaces to point clouds sampled from benchmark B-spline surfaces obtained from structural simulations of heart valve motion. 

These structural simulations are carried out using an IGA-based tricuspid valve parameterization framework \cite{johnson2021param}. This framework defines the annulus during the systolic phase of the cardiac cycle and determines the commissure and leaflet locations along the annular curve. These locations are then used to define the valve heights, whose endpoints are interpolated to construct the free-edge curve of the leaflet. From the annular and free-edge curves, a Gordon surface~\cite{PiegTil97} is generated to reconstruct the leaflets' surface parameterization. For the chordae tendineae, we define the positions of the papillary muscles in $R^{3}$ space, by specifying the number and spacing of the chordae. A steady-state structural simulation is performed on a dense analysis mesh by applying a follower pressure load normal to the leaflet surface, continuing until valve closure. Further details on the surface parameterization and simulation methodology can be found in~\cite{johnson2021param}. The resulting simulation-derived point clouds represent idealized yet exact, noise-free surfaces, where accurate fitting results at different stages of valve motion are expected. To match the point cloud density typically observed in image segmentations of patient data, we sample approximately $1500$ points from the simulation surface for fitting. The sampling is performed using the point cloud simplification filter in MeshLab based on the Poisson disk strategy \cite{cignoni2008meshlab, corsini2012efficient}. 

\begin{figure}[htb!]
    \centering   \includegraphics[width=\linewidth]{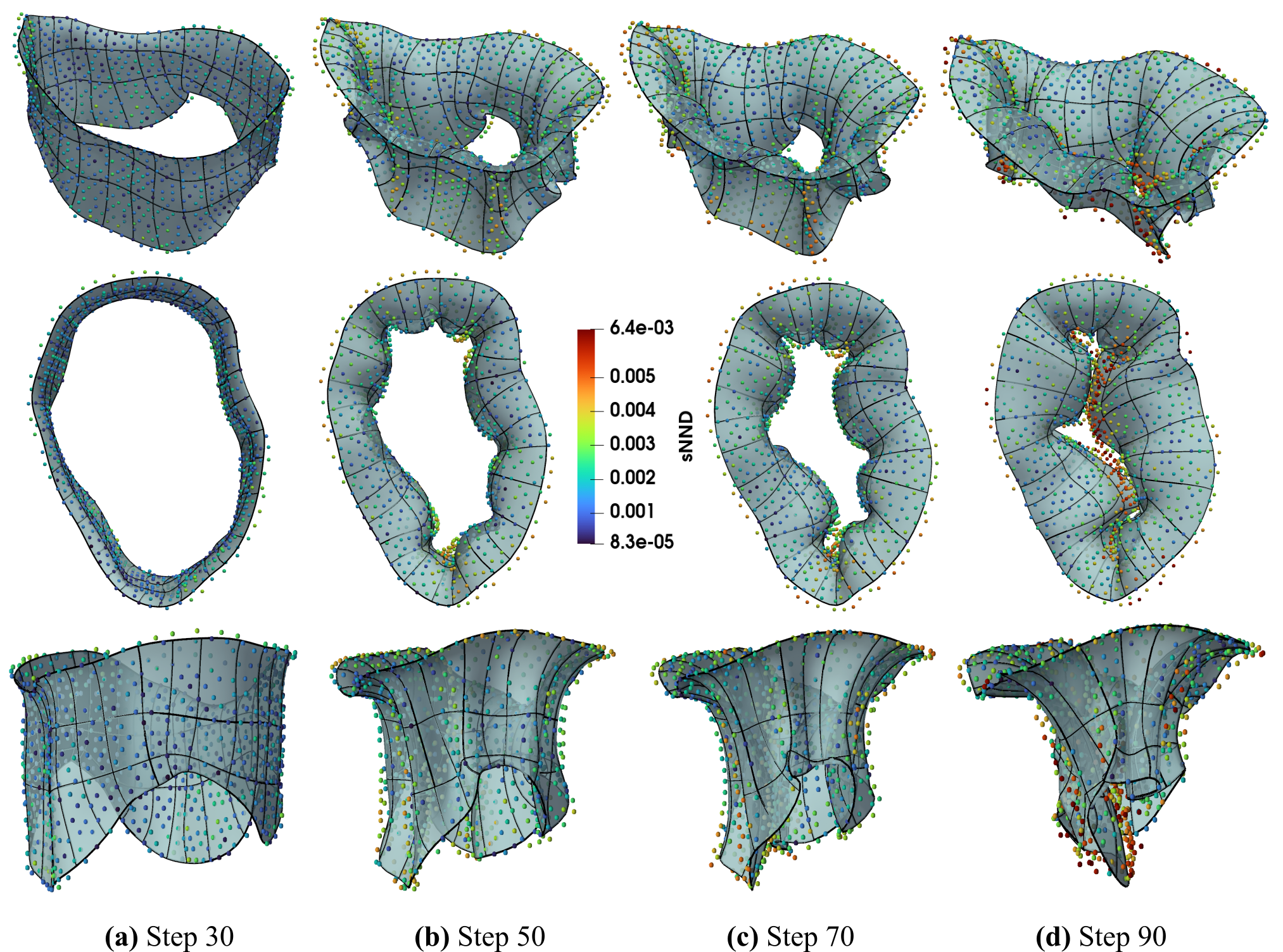}
    \caption{The fitted B-spline surfaces are shown overlaid with the corresponding benchmark point clouds, obtained at multiple steps during the structural simulation of valve motion: \textbf{(a)} Step 30, \textbf{(b)} Step 50, \textbf{(c)} Step 70, and \textbf{(d)} Step 90. The point clouds are colored based on the sNND error with respect to the fitted surface.}
    \label{fig:validation_study}
\end{figure}

The template geometry shown in Fig.~\ref{fig:template} is first re-parameterized using uniform knot vectors, resulting in $33$ control points with degree $p=3$ along the $u$ direction and $6$ control points with degree $q=3$ along the $v$ direction. The template surface is then aligned to the target point cloud at simulation step $30$ (Fig. \ref{fig:validation_study}(a)) using affine transformations. The positions of the control points are then updated by minimizing the loss function, and the resulting fitted surface is used as the template surface for the next point cloud in the cardiac cycle. The weight parameters in the loss function (Eq.~\eqref{eq:loss}) were set to $w_\text{CD} = 80$, $w_\text{HD} = 30$, $w_\text{a} = 0$, $w_\text{orth} = 5$, $w_\text{TPE}=10^{-7}$, and $w_\text{norm} = 2$. These weights are constant during the sequential fitting of the target point clouds. The value of $w_\text{TPE}$ is very low because the TPE term is inherently scale-dependent. In the patient dataset (Section \ref{subsec:patient-data}), the geometries are larger in scale than those in the validation dataset (Section \ref{subsec:validation-study}), resulting in numerically larger inter-point distances. Since TPE values decay sharply with increasing distance (Fig.~\ref{fig:TPE}(b)), we increase the $w_\text{TPE}$ value to maintain an effective repulsive force for larger inter-point distances. The weight $w_\text{a}$ is set to zero because this dataset does not have explicit labels for the points belonging to the annulus and leaflets. We utilize the same weights for the results presented in the Sections~\ref{subsec:density-study}--\ref{subsec:noise-study} as well. 

In Fig.~\ref{fig:validation_study}, we show the fitted B-spline surfaces overlaid with the corresponding point cloud data at each simulation step. The sNND error is computed for each point in the point cloud to evaluate local fitting accuracy. Based on the results, we can see that the fitted surfaces are smooth, intersection-free, and closely aligned with the point clouds across all the stages of valve motion. Complex folding structures of the geometry are effectively captured, while maintaining smoothness and surface regularity (Fig.~\ref{fig:validation_study}(c)-(d)).

\begin{figure}
    \centering
    \includegraphics[width=0.66\linewidth]{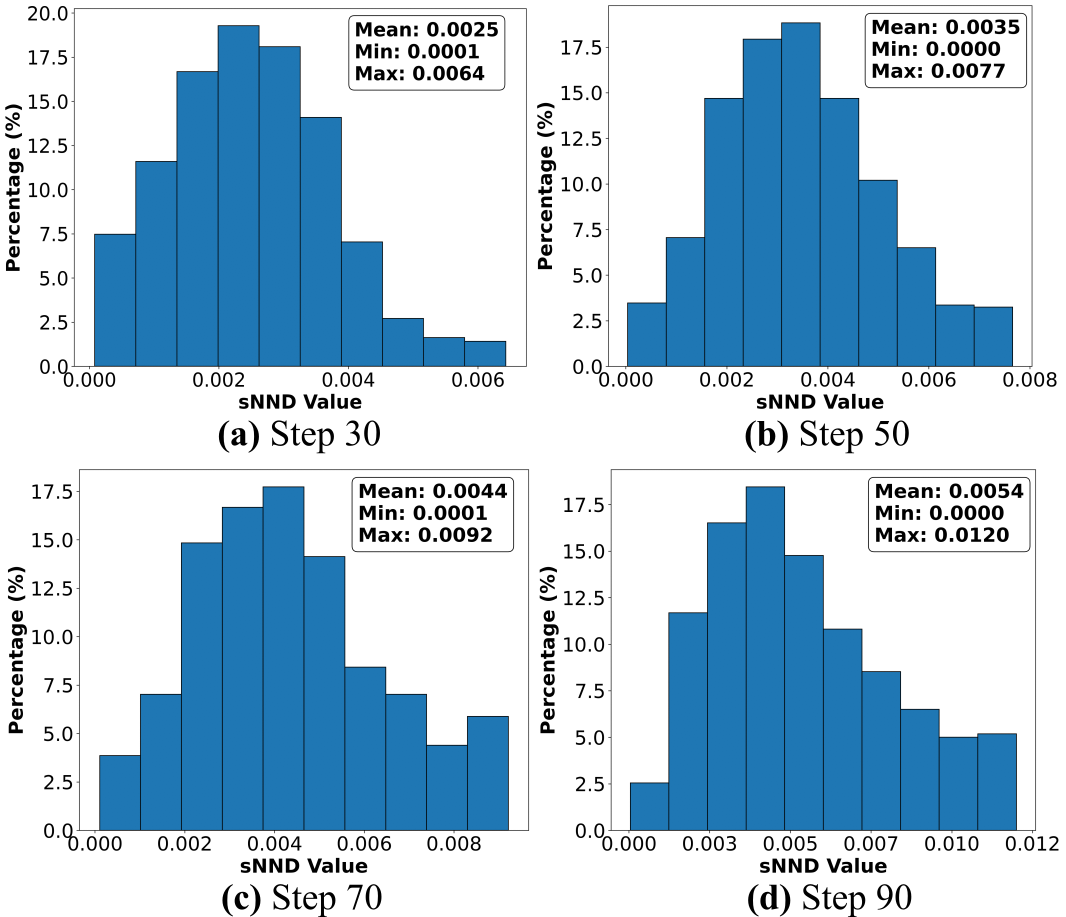}
    \caption{Histograms showing the percentage of points with respect to sNND error at different simulation steps for the validation study. The distributions are consistent across all steps, with the majority of points exhibiting low fitting error, concentrated in the range below $0.005$. This indicates the stable performance of the \textsc{ValveFit} framework over time.}
    \label{fig:validation_histograms}
\end{figure}

\begin{table}[htb!]
    \centering
    \caption{Minimum, maximum, and mean sNND errors between the fitted surface and the target point clouds at different simulation steps.}
    \label{tab:nnd_values_time_series}
    \begin{tabular}{|c|c|c|c|}
        \hline
        \textbf{Step} & \textbf{Minimum sNND} & \textbf{Maximum sNND} & \textbf{Mean sNND} \\
        \hline
        $30$ & $8.25 \times 10^{-5}$ & $6.43 \times 10^{-3}$ & $2.51 \times 10^{-3}$ \\
        $50$ & $3.87 \times 10^{-5}$ & $7.66 \times 10^{-3}$ & $3.47 \times 10^{-3}$ \\
        $70$ & $1.06 \times 10^{-4}$ & $9.21 \times 10^{-3}$ & $4.35 \times 10^{-3}$ \\
        $90$ & $4.37 \times 10^{-5}$ & $1.19 \times 10^{-2}$ & $5.38 \times 10^{-3}$ \\
        \hline
    \end{tabular}
\end{table}

The histograms shown in Fig.~\ref{fig:validation_histograms} show consistent sNND error distributions across all the simulation steps, with the majority of points concentrated in the low sNND range (below $0.005$). Table~\ref{tab:nnd_values_time_series} presents a quantitative assessment of the surface fitting accuracy using the sNND error metric, where lower values indicate better alignment between fitted surfaces and the target point cloud. Here, we report the minimum, mean, and maximum sNND errors for each fitting result. The mean sNND error quantifies the overall fitting accuracy, while the difference between the maximum and minimum sNND errors quantifies the amount of variability in the fitting accuracy. Across all the simulation steps, low sNND error values are observed, with minimum sNND values ranging between $3.87\times10^{-5}$ and $1.06\times10^{-4}$. In all the simulation steps, the maximum sNND values stay below $0.0119$. The mean sNND values range between $2.51\times10^{-3}$ and $5.38\times10^{-3}$, showing consistent fitting at all stages of valve motion. We can see that high accuracy is maintained across most of the geometry, with maximum sNND values localized to regions of greater complexity, particularly near leaflet edges and areas undergoing coaptation as seen in Fig.~\ref{fig:validation_study}(d). The fitted surface also captures the narrow folds seen in Fig.~\ref{fig:validation_study}(c)-(d). The framework effectively balances geometric accuracy with the requirement of maintaining intersection-free surfaces throughout the cardiac cycle.

%%%%%%%%%%%%%%%%%%%%%%%%%%%%%%%%%%%%%%%%%%%%%%%%%%%%%%%%%%%%%%%%%%%%%%%%%%%%%%%%%
\subsection{Surface Fitting Sensitivity with Respect to Point Cloud Sampling}
\label{subsec:density-study}
\begin{figure}
    \centering
    \includegraphics[width=\linewidth]{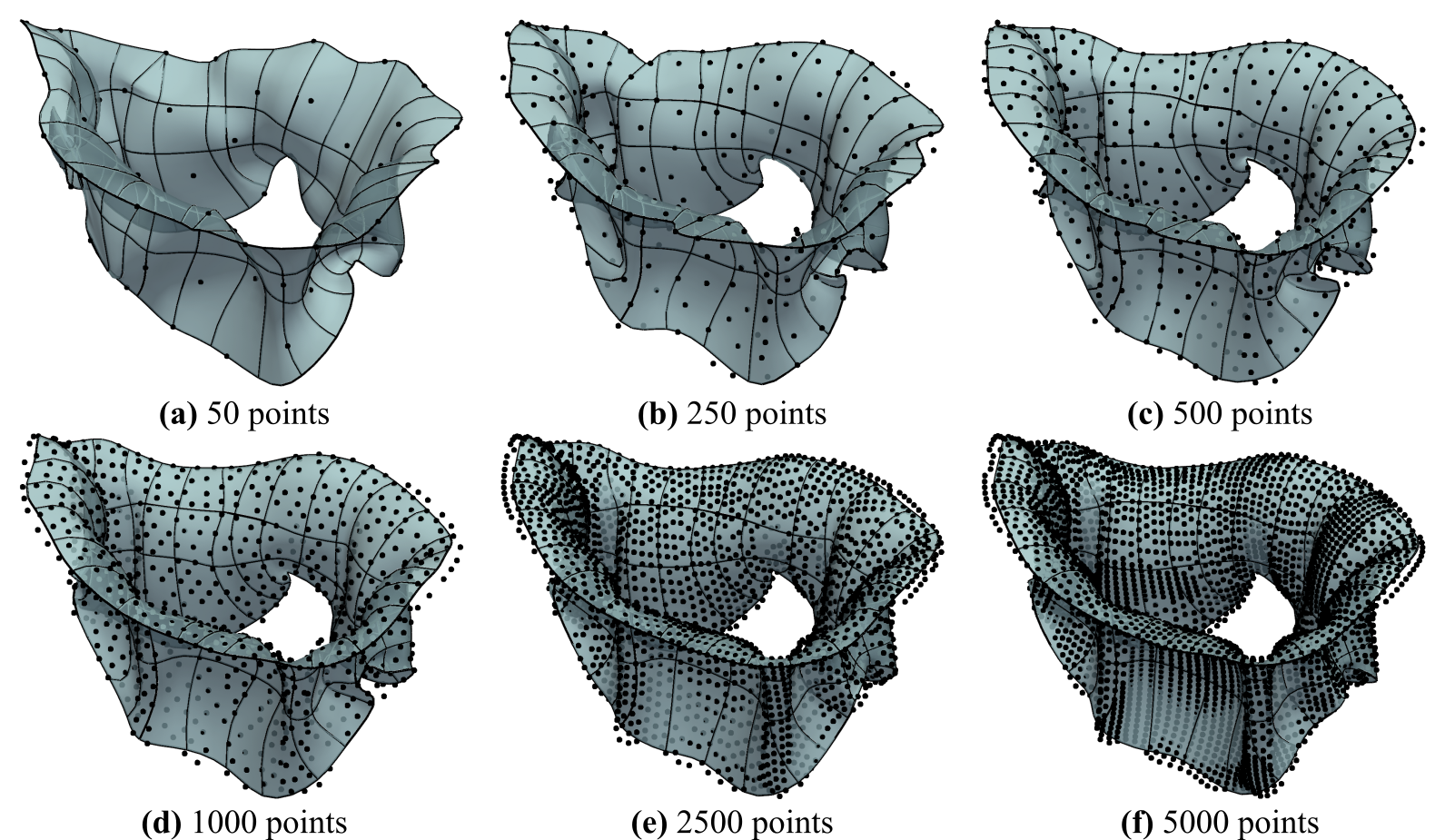}
    \caption{Fitted B-spline surfaces overlaid with target point clouds obtained at simulation step $60$, sampled at varying densities (\textbf{(a)}--\textbf{(f)}: $50$ to $5000$ points).}
    \label{fig:sampling_study}
\end{figure}
We further evaluate the sensitivity of the \textsc{ValveFit} framework to point cloud sampling density by fitting valve surfaces to point clouds of increasing densities. For this study, we performed surface fitting on the point cloud corresponding to the simulation step $60$ of valve motion. The fitted B-spline surfaces are overlaid with the sampled point clouds, as shown in Fig.~\ref{fig:sampling_study}. For each point cloud density, Table~\ref{tab:nnd_values_point_cloud_sampling_updated} reports the minimum, maximum, and mean sNND error values between the fitted surfaces and their corresponding point clouds.

\begin{figure}
    \centering
    \includegraphics[width=\linewidth]{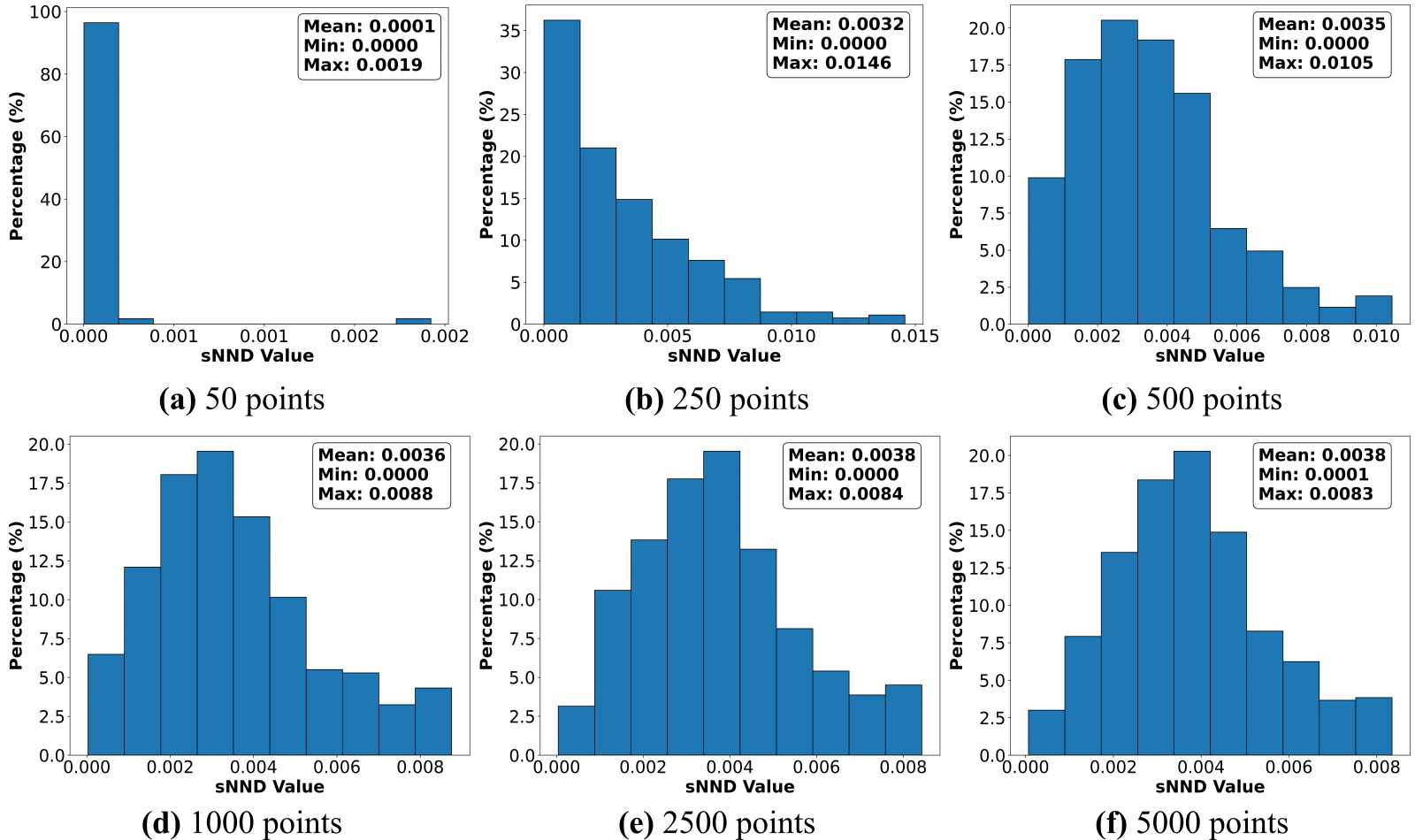}
    \caption{Histograms showing the percentage of points with respect to sNND error at different point cloud densities for simulation data.}
    \label{fig:sampling_study_histograms}
\end{figure}

Based on the results shown in Table~\ref{tab:nnd_values_point_cloud_sampling_updated} and Fig.~\ref{fig:sampling_study_histograms}, we observe that the mean sNND error increases from $6.45\times10^{-5}$ to $0.00376$ as the point cloud density increases from 50 to 2500 points. At a low density of 50 points, the sparse sampling inadequately captures the complexity of the valve geometry, leading to an artificially low sNND error during fitting. As the point cloud density increases, it better represents the intricate geometric features of the valve. Thus, the surface fitting framework needs to resolve more complex geometric features, resulting in a higher and more representative mean error. We note that when calculating the one-sided Chamfer distance---which involves nearest neighbor distances from the point cloud to the surface points---the density of template surface points becomes the limiting factor. Further improvements in the sNND metric due to an increase in point cloud density have no noticeable impact on nearest neighbor distance values once the point cloud density becomes greater than the number of points sampled on the B-spline surface. The results show that the overall topology of the fitted B-spline surface is established even at low point cloud densities. While finer geometric features are resolved as point density increases, the underlying shape of the surface fitted to just 50 points remains  qualitatively consistent with surfaces obtained from denser point clouds.

\begin{table}
    \centering
    \caption{Minimum, maximum, and mean sNND errors between the fitted B-spline surface and target point clouds at varying densities.}
    \label{tab:nnd_values_point_cloud_sampling_updated}
    \begin{tabular}{|c|c|c|c|}
        \hline
        \textbf{Point Count} & \textbf{Minimum sNND} & \textbf{Maximum sNND} & \textbf{Mean sNND} \\
        \hline
        $50$   & $2.55\times 10^{-6}$ & $1.92\times 10^{-3}$ & $6.45\times 10^{-5}$ \\
        $250$  & $5.70\times 10^{-6}$ & $1.46\times 10^{-2}$ & $3.15\times 10^{-3}$ \\
        $500$  & $2.01\times 10^{-5}$ & $1.05\times 10^{-2}$ & $3.50\times 10^{-3}$ \\
        $1000$ & $3.07\times 10^{-5}$ & $8.75\times 10^{-3}$ & $3.57\times 10^{-3}$ \\
        $2500$ & $4.11\times 10^{-5}$ & $8.43\times 10^{-3}$ & $3.76\times 10^{-3}$ \\
        $5000$ & $7.03\times 10^{-5}$ & $8.35\times 10^{-3}$ & $3.82\times 10^{-3}$ \\
        \hline
    \end{tabular}
\end{table}

%%%%%%%%%%%%%%%%%%%%%%%%%%%%%%%%%%%%%%%%%%%%%%%%%%%%%%%%%%%%%%%%%
\subsection{Fitting Sensitivity with Respect to Noise}
\label{subsec:noise-study}
\begin{figure}[H]
    \centering
    \includegraphics[width=\linewidth]{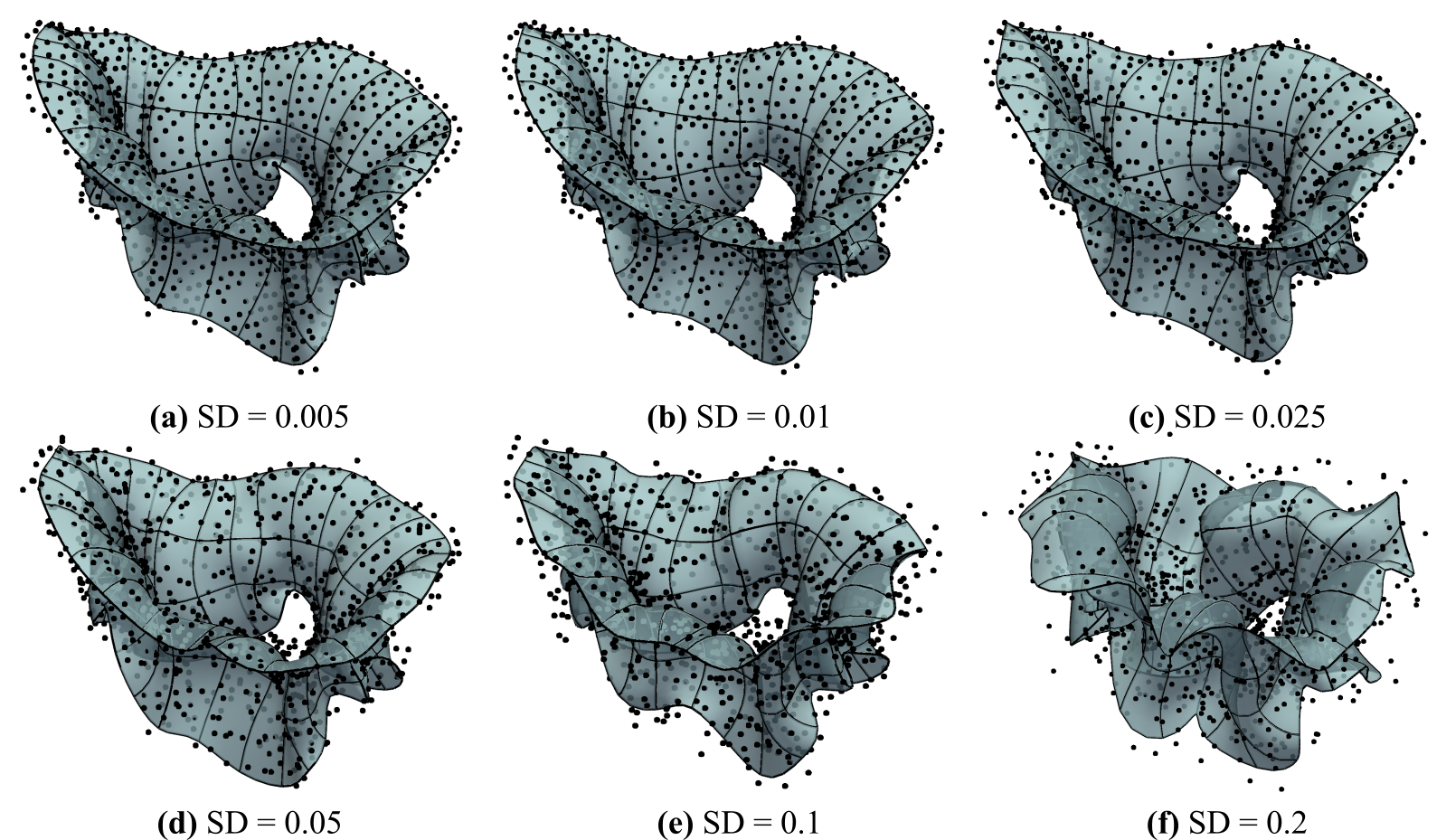}
    \caption{Fitted B-spline surfaces overlaid with target noisy point clouds obtained at simulation step $60$, Gaussian noise at increasing standard deviations is added: \textbf{(a)}--\textbf{(f)} correspond to SD values of $0.005$, $0.01$, $0.02$, $0.05$, and $0.2$, respectively. Despite increasing noise levels, fitted B-spline surfaces are smooth, accurate, and intersection-free.}
    \label{fig:uncolored_noise_study}
\end{figure}

\begin{figure}
    \centering
    \includegraphics[width=\linewidth]{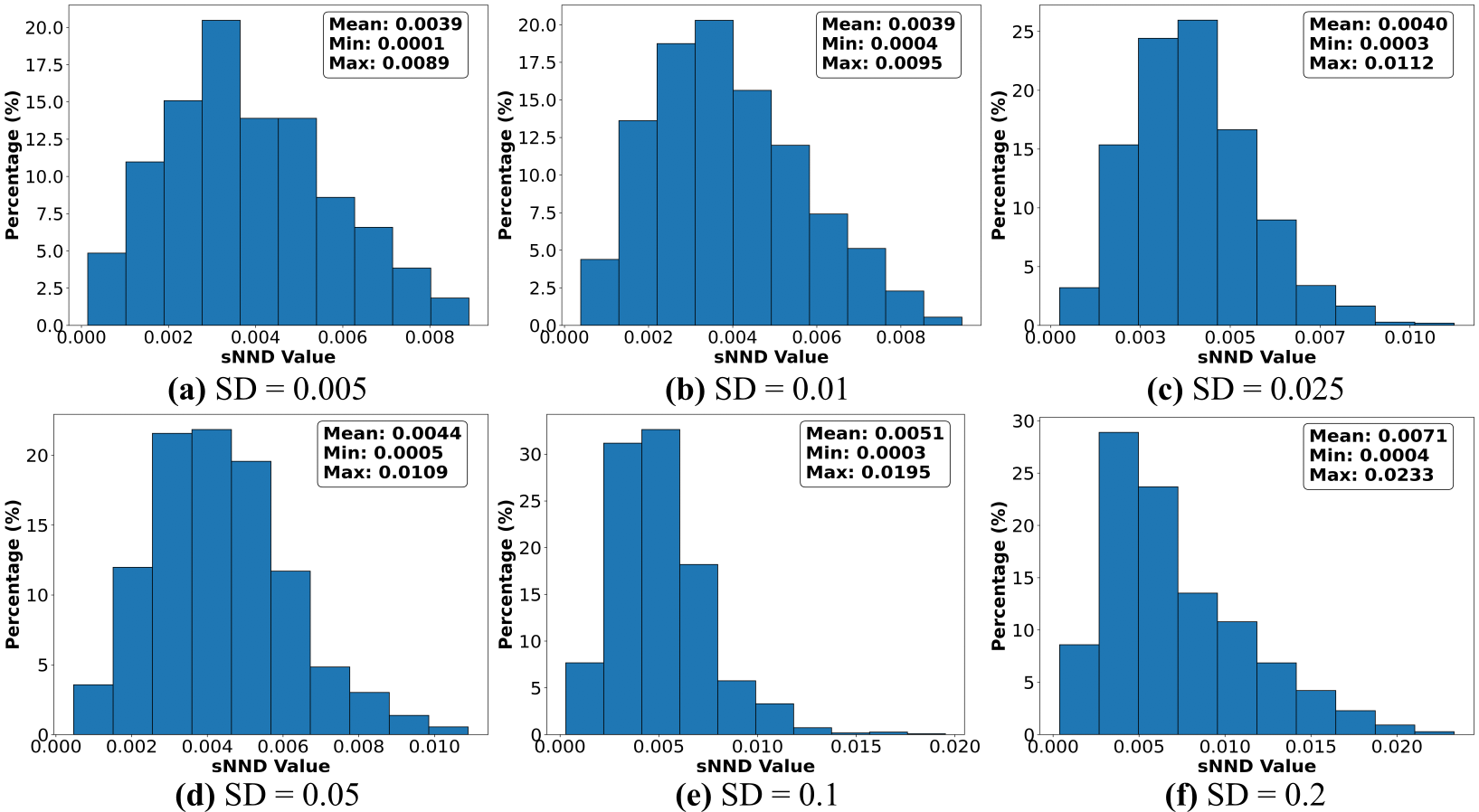}
    \caption{Histogram showing the percentage of points with respect to sNND error at different point cloud noise levels (SD from $0.005$ to $0.2$) for simulation data. We observe that while the spread of errors slightly increases with increasing noise levels, the majority of the points lie in the lower error range.}
    \label{fig:noise_study_histograms}
\end{figure}

In Fig.~\ref{fig:uncolored_noise_study}, we evaluate the robustness of the \textsc{ValveFit} framework under varying levels of noise in the point cloud data. Here, random perturbations were added to the point cloud consisting of $2500$ points at simulation step $60$.  The goal was to simulate the noise typically observed in point cloud data derived from echocardiographic images (Section~\ref{sec:imageacq}). Gaussian noise was added with increasing standard deviations (SD) of $0.005, 0.01, 0.025, 0.05, 0.1, \text{and } 0.2$ \cite{cad2f2ec1bab4deeb7109f8669b42e2b}, representing a specific degree of perturbation. Despite increasing noise, the \textsc{ValveFit} framework successfully captures the finer features of the surface while simultaneously generating smooth and intersection-free surfaces. Even while the noise levels increase, the fitted surfaces maintain a close approximation of the underlying geometry without overfitting. The mean sNND values gradually increase from $0.0039$ to $0.0071$ as the noise levels rise, according to the sNND errors in Table~\ref{tab:noise_impact_snnd_updated}. In Fig.~\ref{fig:noise_study_histograms}, we similarly observe a few points falling into high error bins as the noise levels rise. However, the sNND error does not significantly increase, showcasing minimal distortion of the B-spline surfaces due to noise.

\begin{table}
    \centering
    \caption{Minimum, maximum, and mean sNND errors between the fitted B-spline surface and target point clouds with added Gaussian noise at increasing standard deviations.}
    \label{tab:noise_impact_snnd_updated}
    \begin{tabular}{|l|c|c|c|}
        \hline
        \textbf{SD} & \textbf{Minimum sNND} & \textbf{Maximum sNND} & \textbf{Mean sNND} \\
        \hline
        $0.005$ & $1.48\times 10^{-4}$ & $8.90\times 10^{-3}$ & $3.86\times 10^{-3}$ \\
        $0.01$  & $3.90\times 10^{-4}$ & $9.45\times 10^{-3}$ & $3.90\times 10^{-3}$ \\
        $0.025$ & $2.53\times 10^{-4}$ & $1.12\times 10^{-2}$ & $3.99\times 10^{-3}$ \\
        $0.05$  & $4.67\times 10^{-4}$ & $1.09\times 10^{-2}$ & $4.35\times 10^{-3}$ \\
        $0.1$   & $2.99\times 10^{-4}$ & $1.95\times 10^{-2}$ & $5.10\times 10^{-3}$ \\
        $0.2$   & $3.81\times 10^{-4}$ & $2.33\times 10^{-2}$ & $7.14\times 10^{-3}$ \\
        \hline
    \end{tabular}
\end{table}
In Sections~\ref{subsec:validation-study}--\ref{subsec:noise-study}, we demonstrate the capability of the \textsc{ValveFit} framework to accurately fit point cloud data under different stages of valve motion, point cloud densities, and Gaussian noise levels. When fitting point cloud data sequentially, the framework consistently maintains mean sNND error values between $0.0025$ and $0.0054$, even during complex valve motion. Table~\ref{tab:nnd_values_point_cloud_sampling_updated} shows that the \textsc{ValveFit} framework obtained a mean sNND value of $0.0032$, even with a relatively sparse point cloud density of only $250$ points, which is comparable to the mean sNND values obtained at higher point cloud densities. As can be observed in Table~\ref{tab:noise_impact_snnd_updated}, we also report a gradual increase in sNND values as the amount of Gaussian noise in the point cloud increases. This validation study shows that the \textsc{ValveFit} framework exhibits fitting robustness under a range of conditions by accurately reconstructing surfaces that closely match the ground truth geometry across various noise levels and point cloud densities.

\subsection{Surface Fitting for Patient Data}
\label{subsec:patient-data}
\begin{figure}
    \centering
    \includegraphics[width=0.9\linewidth]{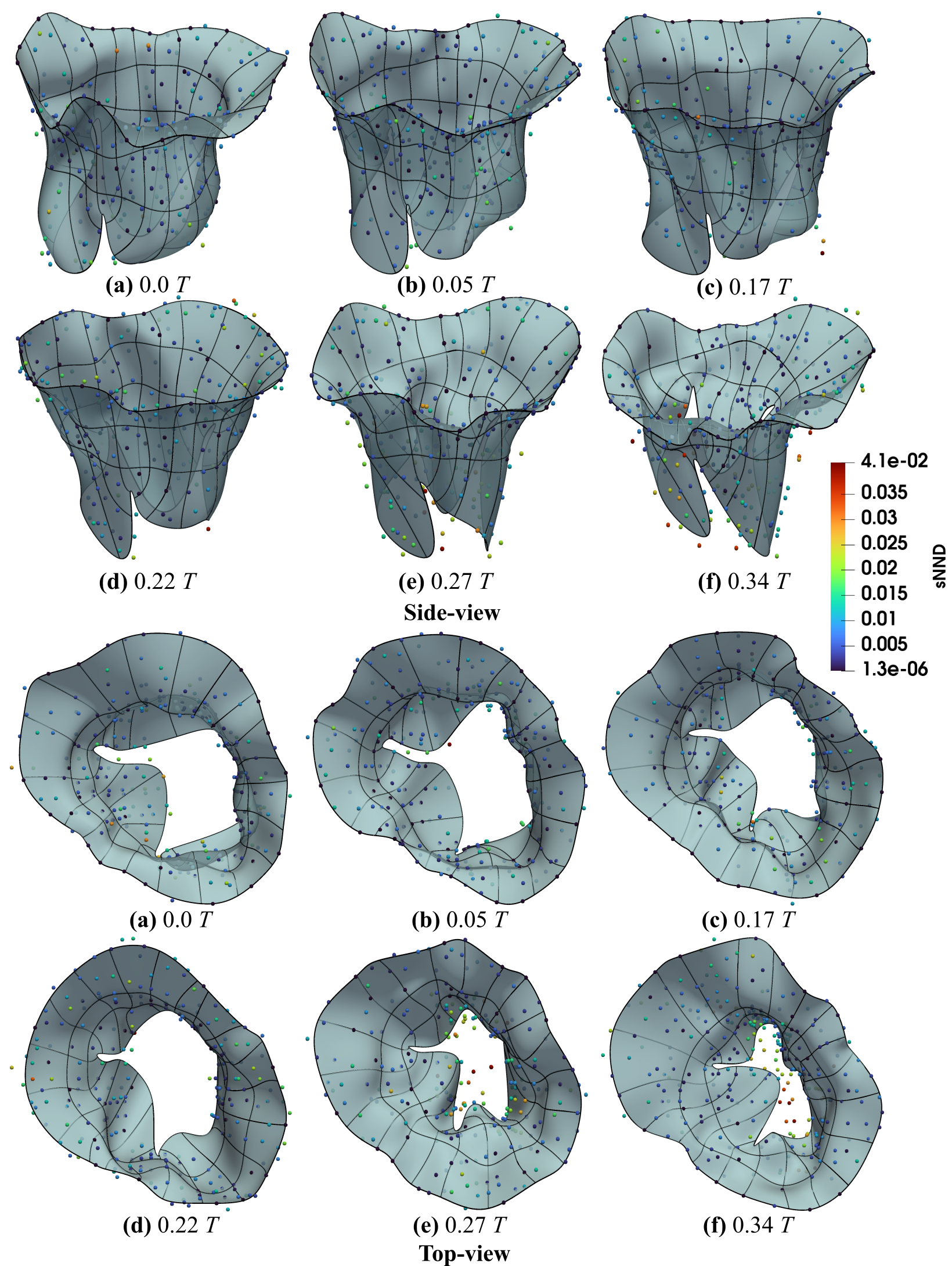}
    \caption{\textbf{Patient 1 dataset:} Fitted surfaces overlaid with the target point clouds at six stages of the cardiac cycle: \textbf{(a)} At the beginning of cardiac cycle $0.0\, T$, \textbf{(b)} $0.05\, T$, \textbf{(c)} $0.17\, T$, \textbf{(d)} $0.22\, T$, \textbf{(e)} $0.27\, T$, and \textbf{(f)} and $0.34\, T$. The first two rows show the side view of the fitted surfaces, while the bottom two rows show the top view of the fitted surfaces. Each point in the target point cloud is colored according to the sNND error.}
    \label{fig:P1_combined_views}
\end{figure}

\begin{figure}
    \centering   \includegraphics[width=0.9\linewidth]{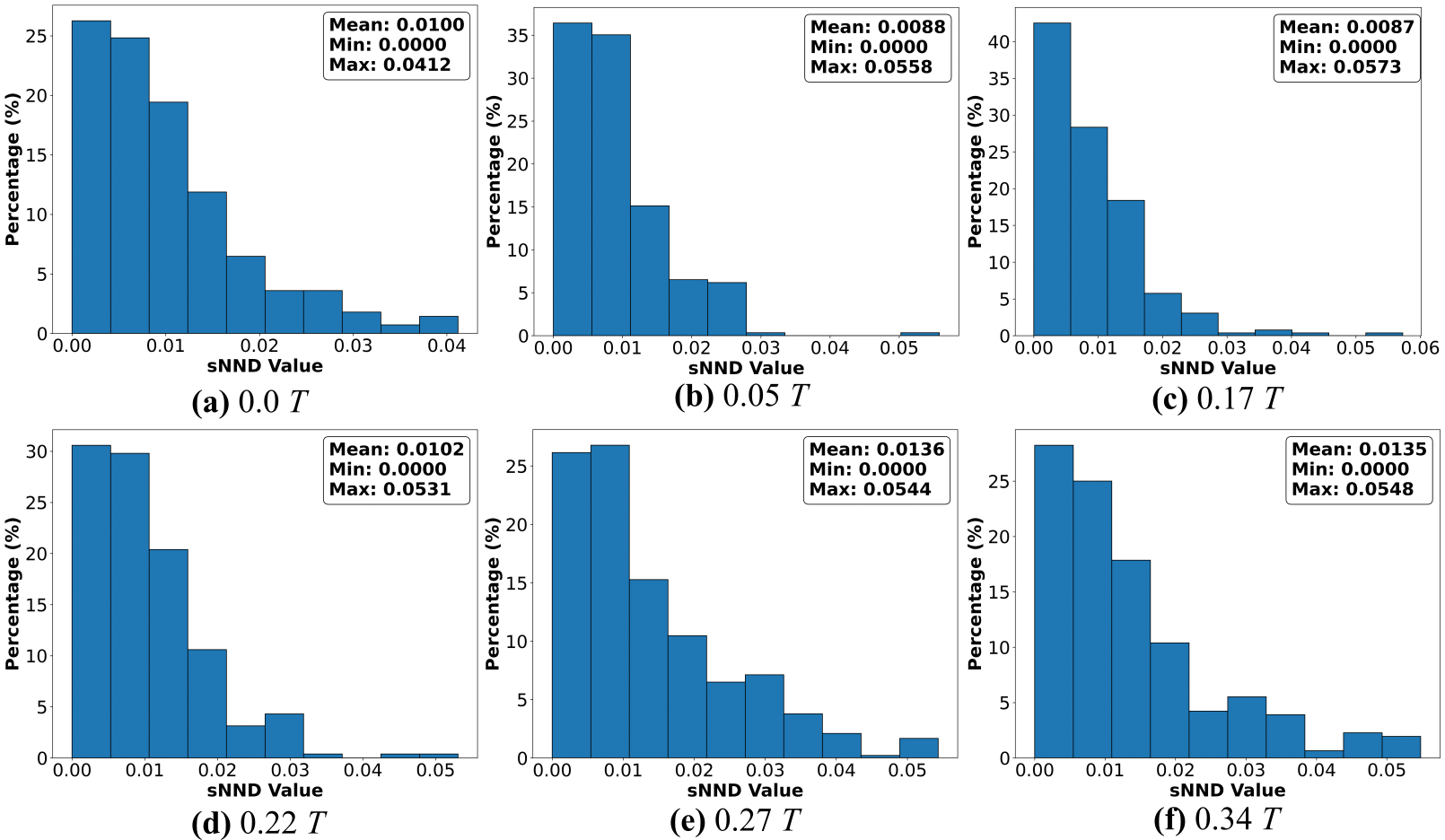}
    \caption{\textbf{Patient 1 dataset:} Percentage of points in the dataset plotted against the sNND error at different time points across the cardiac cycle. Each subplot corresponds to fitting accuracy at a certain time point in the cardiac cycle.}
    \label{fig:P1_histograms}
\end{figure}

\begin{figure}
    \centering
    \includegraphics[width=0.87\linewidth]{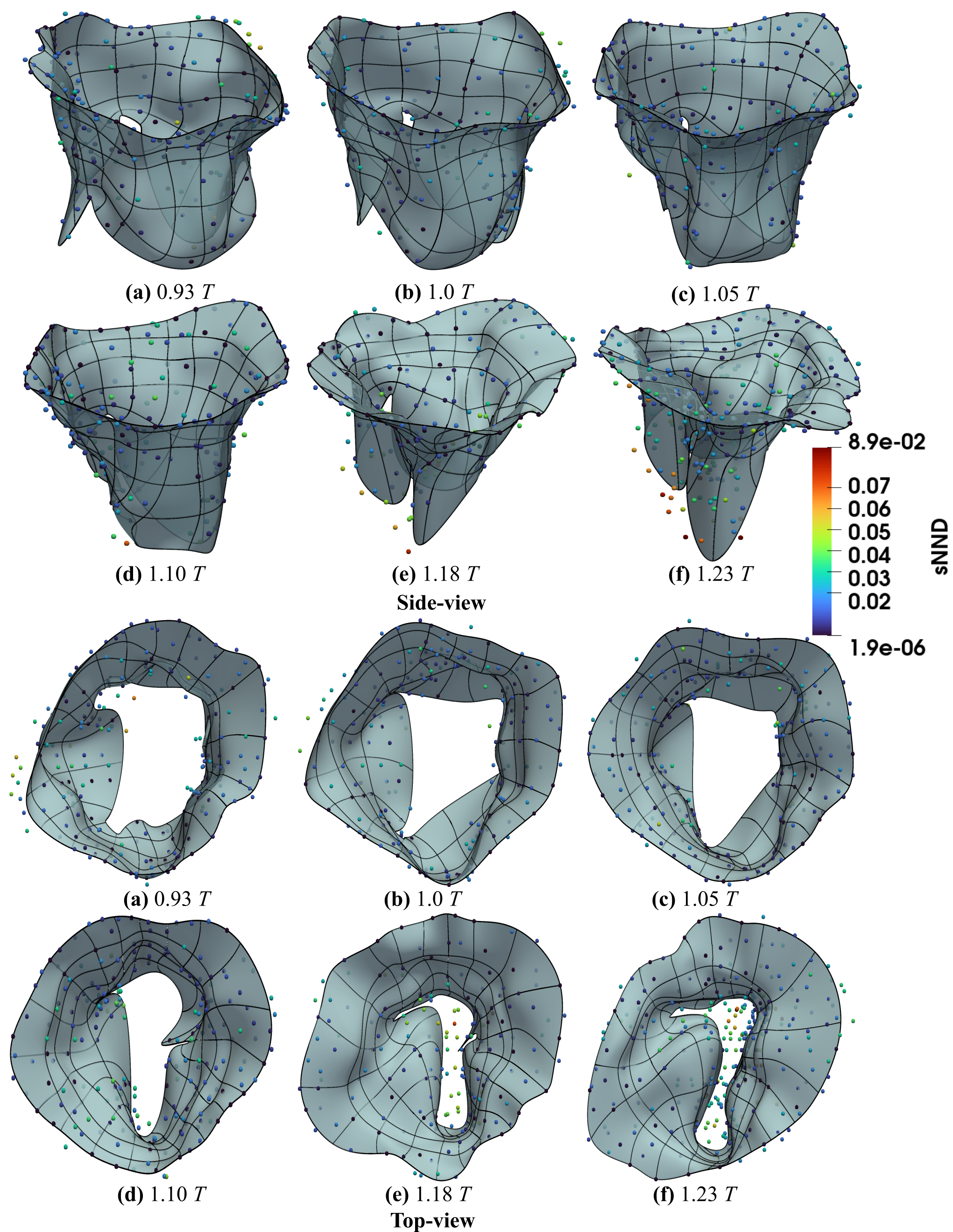}
    \caption{\textbf{Patient 2 dataset:} Fitted TV B-spline surfaces overlaid with the corresponding target point clouds captured at 5 stages of the cardiac cycle: \textbf{(a)} $0.93\, T$, \textbf{(b)} $1.0\, T$, \textbf{(c)} $1.05\, T$, \textbf{(d)} $1.10\, T$, $1.18\, T$, and \textbf{(e)} at the end of the cardiac cycle ($1.23\, T$). The first two rows show the top view of the surfaces, while the last two rows show the side view. Each point in the target point cloud is colored according to its sNND error, providing a visual representation of the accuracy of the fitted surface at that location.}
    \label{fig:P2_combined_views}
\end{figure}

\begin{figure}
    \centering
    \includegraphics[width=\linewidth]{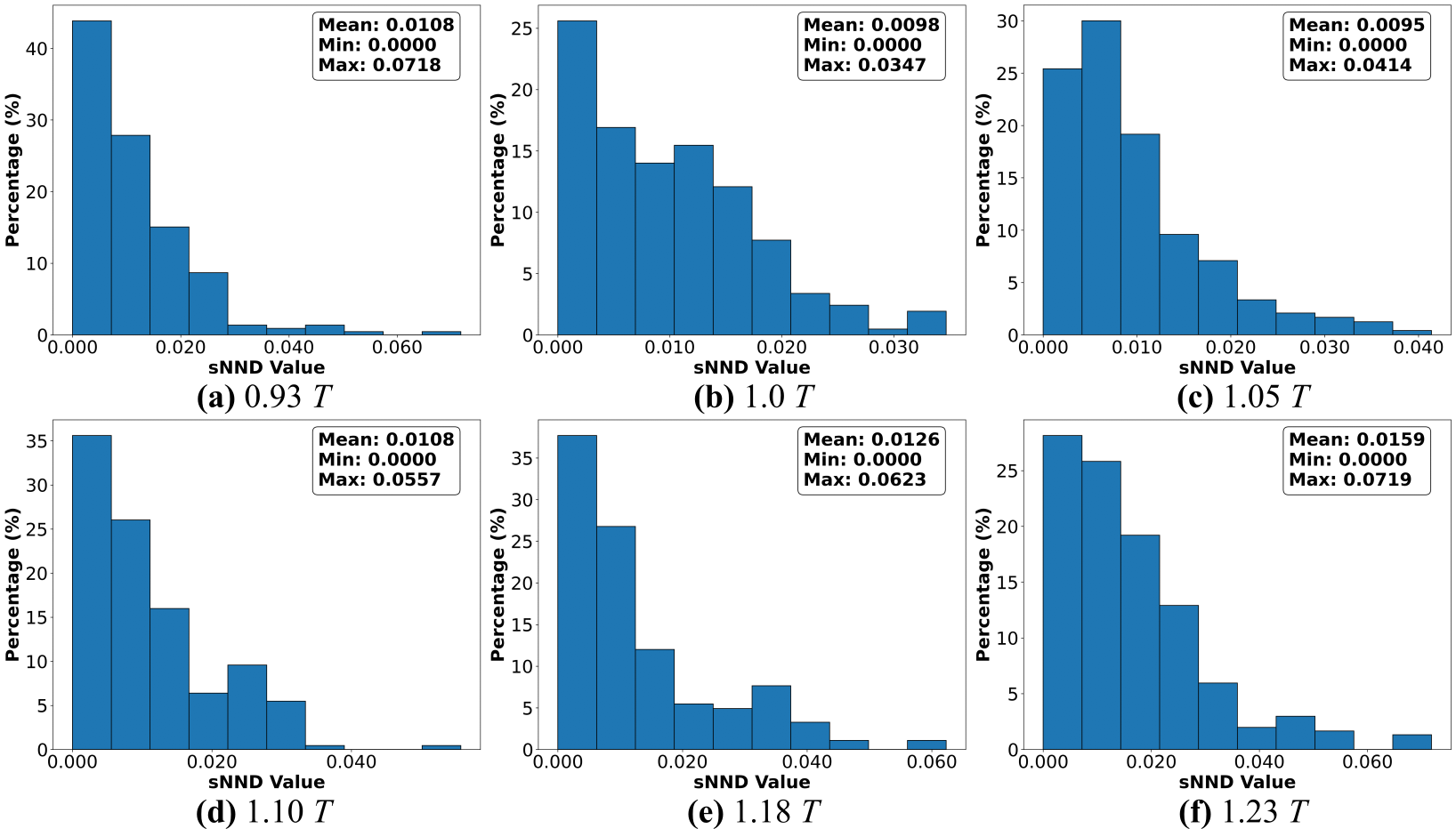}
    \caption{\textbf{Patient 2 dataset:} Percentage of points in the dataset plotted against the sNND error at different time points across the cardiac cycle. Each subplot corresponds to fitting accuracy at a certain time point in the cardiac cycle.}
    \label{fig:P2_histograms}
\end{figure}

After validating the VALVEFIT framework in Sections \ref{subsec:validation-study}-\ref{subsec:noise-study}, we demonstrate its performance for point clouds obtained from echocardiographic imaging data. The point cloud datasets are labeled based on relative time, normalized with the duration of one cardiac cycle ($T$). In this normalized scale, $0.0\, T$ is defined as the beginning of the cardiac cycle, corresponding to the time of minimum right ventricular pressure. $1.0\, T$ represents the end of the same cardiac cycle and corresponds to the beginning of the next cardiac cycle ($1.0\, T_{i} = 0.0\, T_{i+1}$). Important physiological events are also marked within a cardiac cycle, where $0.3\, T$ corresponds to end-diastole (ED), when the valve leaflets are coapted, and approximately at $0.45\, T$, which corresponds to end-isovolumetric contraction \cite{nguyen2019dynamic}. Figs.~\ref{fig:P1_combined_views} and \ref{fig:P1_histograms} show surface fitting results for the first patient data at $0.0\, T$, $0.05\, T$, $0.17\, T$, $0.22\, T$, $0.27\, T$, and $0.34\, T$, while Figs.~\ref{fig:P2_combined_views} and \ref{fig:P2_histograms} show the results for the second patient data at $0.93\, T$, $1.0\, T$, $1.05\, T$, $1.10\, T$, $1.18\, T$, and $1.23\, T$.

The template geometry shown in Fig.~\ref{fig:template} is used for fitting patient data. The template surface parameterization for fitting the Patient 1 dataset consists of $23$ control points with degree $p=3$ basis functions along the $u$ direction and $6$ control points with degree $q=3$ basis functions along the $v$ direction, resulting in a control grid of size $23\times6$. The weight parameters in Eq.~\eqref{eq:data-fidelity} for the first patient data were set to $w_\text{CD} = 20$, $w_\text{HD} = 0.5$, $w_\text{a} = 1$, $w_\text{orth} = 10$, $w_\text{TPE}=1.5$, and $w_\text{norm} = 5$. The weights are identical to those used in Section \ref{subsec:validation-study}, except for an increased Tangent Point Energy weight ($w_\text{TPE}$) to account for larger inter-point distances. Additionally, for patient data, $w_\text{a}$ is set to a non-zero value because we have access to point clouds where the annulus and leaflets have been separately segmented. For the Patient 2 dataset, we observed greater geometric complexity and higher noise intensity in the point clouds. The template surface parameterization is changed to result in a control grid of size $23\times8$. To avoid overfitting to the noise, we thus increased the relative importance of the regularization terms in the loss function (Eq.~\eqref{eq:loss}) by increasing the weights $w_\text{orth}$ to $15$ and $w_\text{norm}$ to $7.5$, ensuring that the surface does not attempt to overfit the noisy point cloud and remains smooth.

\begin{table}[!htb]
    \centering
    \caption{\textbf{Patient 1 dataset:} Minimum, maximum and mean sNND errors between the fitted B-spline surface and target point clouds at different timesteps.}
    \label{tab:P1}
    \begin{tabular}{|c|c|c|c|}
        \hline
        \textbf{Step} & \textbf{Minimum sNND} & \textbf{Maximum sNND} & \textbf{Mean sNND} \\
        \hline
        $0.0\, T$ & $1.2566\times 10^{-6}$ & 0.04120 & 0.00995 \\
        $0.05\, T$ & $1.2560\times 10^{-6}$ & 0.05581 & 0.00881 \\
        $0.17\, T$ & $2.8330\times 10^{-6}$ & 0.05728 & 0.00871 \\
        $0.22\, T$ & $3.0224\times 10^{-6}$ & 0.05310 & 0.01019 \\
        $0.27\, T$ & $1.6789\times 10^{-6}$ & 0.05436 & 0.01357 \\
        $0.34\, T$ & $1.0263\times 10^{-6}$ & 0.05481 & 0.01354 \\
        \hline
    \end{tabular}
\end{table}

In Figs.~\ref{fig:P1_histograms} and \ref{fig:P2_histograms}, we show histograms of the sNND error distributions from the target point clouds to the fitted surface for Patients 1 and 2, respectively. For the Patient 1 dataset (Fig.~\ref{fig:P1_histograms}), the maximum observed sNND error is $0.05728$, with the majority of points reporting sNND errors below $0.03$. The maximum observed sNND error for the Patient 2 dataset (Fig.~\ref{fig:P2_histograms}) is $0.07187$, with a greater percentage of points reporting an sNND error below $0.03$. The distributions show consistent performance of the framework throughout the cardiac cycle. According to Table~\ref{tab:P1}, the mean sNND errors for the Patient 1 dataset are between $0.00995$ and $0.01357$ across the six stages, showing the framework's consistently stable performance throughout the cardiac cycle. The maximum sNND values ranged from $0.0412$ to $0.05728$, with slightly higher values at time points $0.17\, T$ and $0.34\, T$. These values were found to correlate with increased geometric complexity and leaflet coaptation, respectively.

\begin{table}[!htb]
    \centering
    \caption{\textbf{Patient 2 dataset:} Minimum, maximum and mean sNND errors between the fitted B-spline surface and target point clouds at different timesteps.}
    \label{tab:P2}
    \begin{tabular}{|c|c|c|c|}
        \hline
        \textbf{Timestep} & \textbf{Minimum sNND} & \textbf{Maximum sNND} & \textbf{Mean sNND} \\
        \hline
        $0.93 \, T$ & $4.1777 \times 10^{-6}$ & 0.07178 & 0.01082 \\
        $1.0 \, T$ & $1.9923 \times 10^{-6}$ & 0.03468 & 0.00976 \\
        $1.05 \, T$ & $2.1534 \times 10^{-6}$ & 0.04140 & 0.00947 \\
        $1.10 \, T$ & $4.5950 \times 10^{-6}$ & 0.05571 & 0.01082 \\
        $1.18 \, T$ & $2.2556 \times 10^{-6}$ & 0.06232 & 0.01262 \\
        $1.23 \, T$ & $1.1482 \times 10^{-6}$ & 0.07187 & 0.01593 \\
        \hline
    \end{tabular}
\end{table}

\begin{figure}[!htb]
    \centering
    \includegraphics[width=0.95\linewidth]{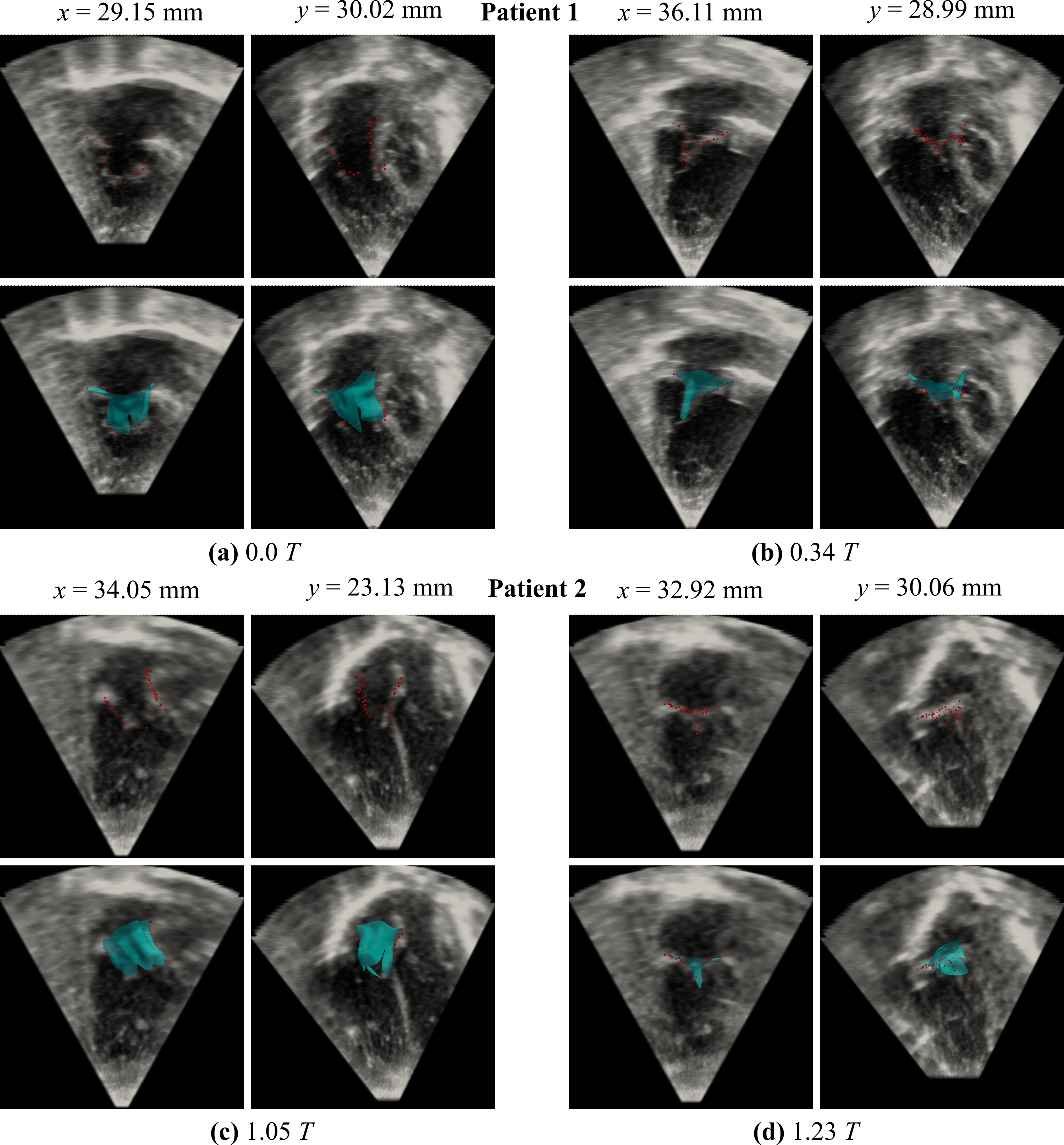}
    \caption{Overlays of the fitted B-spline surfaces (blue) and segmented target point clouds (red) on the corresponding 2D echocardiography slices for Patient~1 at \textbf{(a)}~0.0~T and \textbf{(b)}~0.32~T, and for Patient~2 at \textbf{(c)}~1.05~T and \textbf{(d)}~1.32~T. Each panel demonstrates the agreement of the fitted patient-specific surfaces with the underlying image data.}
    \label{fig:echo_overlay}
\end{figure}

Compared to the Patient 1 dataset, the sNND errors for the Patient 2 dataset (Table~\ref{tab:P2}) were slightly higher, with mean sNND values ranging from $0.00947$ to $0.01593$. The fitting accuracy varies throughout the cardiac cycle. Both the mean and maximum sNND errors decrease initially, reaching a minimum value near the fully open valve configuration around $1.0\,T$ and $1.05\,T$. Subsequently, the errors trend upward through the remainder of the cycle. The Patient 2 dataset required weight parameter adjustments, as discussed above, due to higher noise levels in the point clouds. The reduced relative importance of the Chamfer distance loss term compared to the regularization terms is the primary contributor to the slight increase in sNND error. Through hyperparameter tuning, the \textsc{ValveFit} framework balances geometric regularity and fitting accuracy while dealing with different noise levels. Thus, this framework can be used to generate analysis-suitable patient-specific B-spline surfaces for a variety of patient datasets with varying image quality and segmentation accuracy. In Fig.~\ref{fig:echo_overlay}, we present the fitted B-spline surfaces overlaid on echocardiographic slices, together with the target points segmented from the corresponding slices for the representative time points of Patients~1 and~2 datasets. The slices are extracted along planes oriented orthogonally to either the x-axis (sagittal) or y-axis (coronal), with the slicing positions specified by the coordinates shown above each panel. The overlays demonstrate a close correspondence between the fitted surfaces and the key anatomical structures in the echocardiographic slices, providing visual confirmation of the framework’s accuracy. These results highlight the capability of the proposed approach to generate patient-specific, analysis-suitable surface representations that are consistent with clinical imaging data.

\subsection{Tricuspid Valve Structural Simulation}
\label{sec:simulation}

Accurate modeling of TV closure under systolic loads is vital for predicting regurgitation and supporting surgical planning. To demonstrate the utility of our framework for generating analysis-suitable B-spline surfaces in patient-specific applications, we perform a simplified structural simulation on the fitted B-spline surface and highlight its ability to extract meaningful metrics, such as strain. For the atrioventricular heart valves, the primary driving force for valve opening and closure is the pressure gradient across the valve leaflet surfaces, which can be idealized or realized by means of the pressure loading in the structural mechanics analysis \cite{morganti2015patient, kamensky2018contact, johnson2021param, zhang2021simulating, kiendl2015isogeometric, sun2014computational, toma2016high, sacks2019simulation,votta2008mitral,laurence2020pilot, stevanella2010finite,kong2018finite,mathur2022texas}. In this study, we adopt the isogeometric analysis framework in ~\cite{johnson2021param} to simulate valve closure using a patient-specific surface fitted with \textsc{ValveFit} framework.

For this simulation, the point cloud for Patient 2 at $1.0\, T$ is considered the initial, stress-free configuration.  The leaflet material behavior is modeled using the isotropic Lee–Sacks model~\cite{wu2018a}, where the strain energy functional is defined as
$
    \frac{1}{2} c_0 ( I_1 - 3 ) + \frac{1}{2} c_1 (e^{c_2 (I_1-3)^2} - 1 )
$.
Here, $I_1$ is the first invariant (trace) of the left Cauchy--Green deformation tensor, and $c_0$, $c_1$, and $c_2$ are the material coefficients set to $c_0$ = $10$ kPa, $c_1$ = $0.209$ kPa, and $c_2$ = $9.046$ \cite{kamensky2018a}. The leaflet's mass density is set at $1.0$ g/cm$^3$ with a thickness of $0.0396$ cm.

\begin{figure}
    \centering
    \includegraphics[width=1.0\textwidth]{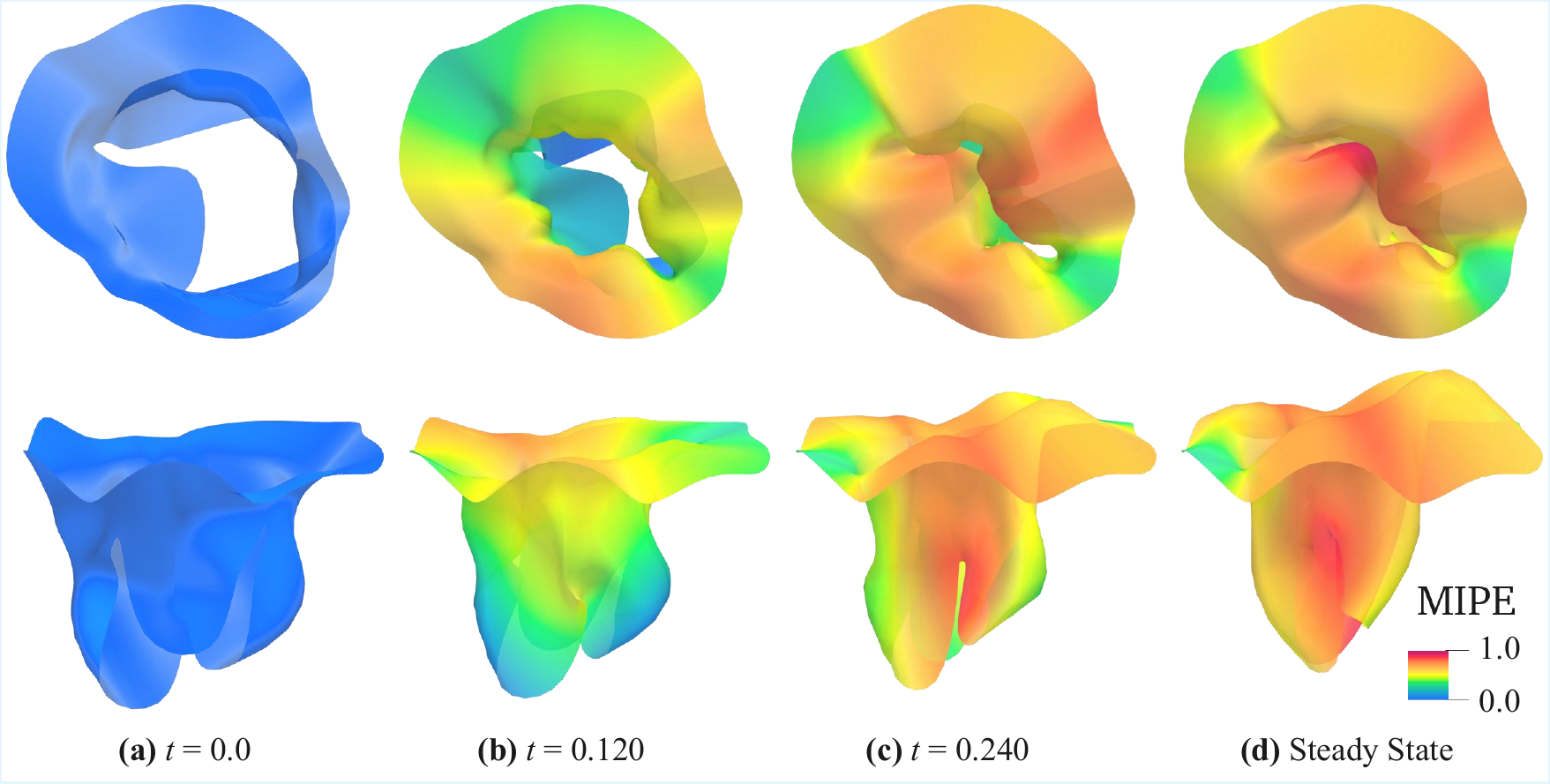}
    \caption{A series of snapshots illustrating the closing behavior of the tricuspid valve from top and side views. The color indicates the maximum in-plane principal Green--Lagrange strain (MIPE) on the atrial side of the leaflets. A pressure load of 25 mmHg increases linearly over the course of 5 milliseconds, mirroring the physiological pressure gradient across the TV during systole for this simulation.}
    \label{fig:TV_Simulation}
\end{figure}

Accurate TV dynamics require modeling the chordae tendineae, as healthy chordae play an important role in preventing valve prolapse. Without them, the valve cannot properly seal, leading to regurgitation and other functional impairments. Despite their importance to valve function, it is difficult for current non-invasive imaging techniques to capture chordae with enough fidelity for accurate geometric modeling. Previous efforts have worked around this issue by creating functional chordae maps~\cite{Narang2021} or by using interactive dynamic simulations to estimate chordae locations heuristically to match patient-valve behavior~\cite{Walczak2021}.

This work uses simplified chordae to prevent leaflet prolapse and demonstrates how the geometry is directly analyzable through closure mechanics. Six idealized cables, each 1 cm in length, are attached to the leaflet and extend outward, perpendicular to its surface. Each TV leaflet has two chordae, with one end of each chordae attached to an insertion point along the free edge of the TV, and the other fixed in space to represent papillary muscle attachment. An IGA cable formulation \cite{johnson2021param} is used for the chordae, which are modeled as cubic B-spline curves with a Young’s modulus of $40$~MPa and a fictitious radius of $0.023$ cm.  This setup ensures that the valve is constrained and resists leaflet prolapse. For more accurate patient-specific simulations, proper insertion points and papillary muscle locations can be used to achieve more clinically viable results.

The simulations specifically focus on the tricuspid valve in its closed state. For simplification, the positions of the annulus and papillary muscles are fixed in the initial configuration of $1.0\, T$, and thus their motion or deformation is not modeled. Valve closure is achieved by applying a pressure load on the ventricular side of the leaflet surface. This load increases linearly from $0$ to $25$ mmHg over the course of $5$ milliseconds ($0$ to $0.005$ s), mirroring the physiological pressure gradients across the TV during systole.

Fig.~\ref {fig:TV_Simulation} shows the top and side views of the simulated valve closure, highlighting the framework's ability to produce patient-specific, analysis-suitable geometries for predictive modeling. These geometries allow for the recovery of valuable cardiac cycle metrics, including deformation, coaptation, and strain (e.g. MIPE, the maximum in-plane principal Green--Lagrange strain), which are often inaccessible through clinical imaging alone. With access to more complete patient data, it should be possible to recover physiologically relevant, patient-specific metrics throughout the entire cardiac cycle.

\section{Conclusion}\label{sec:conclusion}
In this article, we introduce \textsc{ValveFit}, a novel GPU-accelerated B-spline surface fitting framework for 3D surface reconstruction of tricuspid heart valve geometry from noisy and sparse point cloud data obtained via medical image segmentation. In this approach, surface fitting is carried out by deforming a template B-spline geometry. We introduce novel regularization terms that maintain the local orthogonality of the surface tangent vectors, prevent global self-intersections between surfaces through tangent point energy, and constrain the variation of the normal vectors to ensure that the underlying B-spline surface is smooth and intersection-free. The resulting B-spline surfaces are suitable for biomechanical simulations. We first validated the \textsc{ValveFit} framework with a benchmark point cloud dataset that is derived from structural simulations. Here, we also demonstrate robustness toward fitting varying point cloud densities and increasing noise intensities. Finally, we have demonstrated surface fitting for real patient data, where we accurately captured the complex three-leaflet structure of the tricuspid valve. The fitted surfaces closely approximate the inherent valve geometry represented by the point clouds across multiple time points in the cardiac cycle. We have demonstrated the robustness of the proposed ValveFit framework for point clouds obtained from echocardiographic images of heart valves, which are inherently noisy and have a lower spatial resolution compared to other imaging modalities, such as cine MR or CT. Thus, we envision that the ValveFit framework can be easily extended to point cloud data obtained from other imaging modalities as well. We demonstrate that the fitted surfaces are suitable for biomechanical analysis by performing a structural simulation of the fitted geometry and obtaining the strain distribution during valve deformation. These patient-specific valve geometries can offer insight into the progression of valve-related complications, enabling biomechanical modeling-guided precise and effective interventions. 

In our future work, we would like to extend the applicability of the ValveFit framework to patient-specific geometric modeling of the mitral, aortic, and pulmonary heart valves. To apply the proposed framework to these valve types, one simply needs to start from their respective idealized template geometries. Similar to the approach discussed above to construct the template geometry of the tricuspid valve, we can generate anatomical and topologically similar idealized template geometries for the rest of the heart valves, which would serve as the input to the ValveFit framework for patient-specific geometric modeling. We also plan to extend the \textsc{ValveFit} framework to support volumetric spline parameterization and patient-specific volumetric modeling for isogeometric analysis.

\section*{Code Availability}
The code is available at \url{https://github.com/CMIA-Lab/ValveFit.git} under the MIT license. 

\section*{Acknowledgments}
A.~Moola and A.~Pawar would like to acknowledge funding from the Translational AI Center at Iowa State University. A.M.~Corpuz and M.-C.~Hsu would like to acknowledge funding from the National Science Foundation under award number DMS-2436623. C.-H.~Lee would like to acknowledge support from the National Institutes of Health (R01HL159475, R01GM157589), the American Heart Association (16SDG27760143, 24IVPHA1283213), and the Oklahoma Center for the Advancement of Science and Technology (OCAST HR23-03). This support is gratefully acknowledged.

\bibliographystyle{unsrtnat}

\bibliography{ref-mch-r3}

\end{document}